\documentclass{amsproc}
\usepackage{amssymb}
\usepackage{amsmath}
\usepackage[matrix,arrow,curve,cmtip]{xy}

\numberwithin{equation}{subsection}

\theoremstyle{plain}   
\newtheorem{bigthm}{Theorem}   

\newtheorem{cor}[equation]{Corollary}     
\newtheorem{lemma}[equation]{Lemma}         
\newtheorem{prop}[equation]{Proposition} 
\newtheorem{addendum}[equation]{Addendum}
    
\newtheorem*{unconj}{Conjecture}

\theoremstyle{definition}

\theoremstyle{remark}
\newtheorem{remark}[equation]{Remark}

\newcommand{\Fil}{\operatorname{Fil}}
\newcommand{\Tor}{\operatorname{Tor}}

\newcommand{\Hom}{\operatorname{Hom}}

\newcommand{\Ext}{\operatorname{Ext}}

\newcommand{\f}{\hat{\phantom{:}}}
\newcommand{\cf}{\check{\phantom{i}}}
\newcommand{\cs}{^{\boldsymbol{\cdot}}}
\newcommand{\gr}{\operatorname{gr}}
\newcommand{\TF}{\operatorname{TF}}
\newcommand{\TC}{\operatorname{TC}}
\newcommand{\TR}{\operatorname{TR}}
\renewcommand{\TH}{\operatorname{TH}}

\newcommand{\N}{\mathbb{N}}
\newcommand{\Z}{\mathbb{Z}}

\newcommand{\C}{\mathbb{C}}
\newcommand{\Zp}{\mathbb{Z}_p}
\newcommand{\Qp}{\mathbb{Q}_p}
\newcommand{\Fp}{\mathbb{F}_p}
\newcommand{\cont}{{\text{\rm cont}}}

\newcommand{\colim}{\operatornamewithlimits{colim}}
\newcommand{\holim}{\operatornamewithlimits{holim}}
\newcommand{\hocolim}{\operatornamewithlimits{hocolim}}

\newcommand{\id}{\operatorname{id}}
\newcommand{\pr}{\operatorname{pr}}
\newcommand{\tr}{\operatorname{tr}}
\newcommand{\xto}{\xrightarrow}
\newcommand{\e}{\varepsilon}
\newcommand{\tate}{\hat{\mathbb{H}}}
\newcommand{\borel}{\mathbb{H}\hskip1pt_{\boldsymbol{\cdot}}}
\newcommand{\coborel}{\displaystyle{\mathbb{H}\,{}^{\boldsymbol{\cdot}}}}

\begin{document}

\title{On the topological cyclic homology of the algebraic\\ closure of
  a local field}

\author{Lars Hesselholt*}

\address{Massachusetts Institute of Technology, Cambridge,
Massachusetts}

\email{larsh@math.mit.edu}

\address{Nagoya University, Nagoya, Japan}

\email{larsh@math.nagoya-u.ac.jp}

\thanks{${}^*)$\ Supported in part by the National Science
Foundation (USA) and COE (Japan)} 

\maketitle

\section*{Introduction}

The cyclotomic trace map provides a comparison of the algebraic
$K$-theory spectrum with a pro-spectrum $\{\TR^n\}_{n \geq 1}$ that is
built from the cyclic fixed points of topological Hochschild
homology. In the long paper~\cite{hm2} we used this comparison and
an approximate evaluation of the structure of the pro-spectrum
$\{\TR^n\}_{n \geq 1}$ to completely determine the $p$-adic $K$-theory
of a mixed characteristic local field $K$. This verified the
Lichtenbaum-Quillen conjecture for the field $K$. In this paper we
completely determine the structure of the pro-spectrum $\{\TR^n\}_{n
\geq 1}$ for the algebraic closure $\bar{K}$ of the local field
$K$. This leads us to formulate a conjecture for the structure of the
pro-spectrum $\{\TR^n\}_{n \geq 1}$ for the local field $K$.

Let $V$ be a complete discrete valuation ring with quotient field $K$
of characteristic $0$ and with perfect residue field $k$ of odd
characteristic $p$. In the paper~\cite{hm2} we constructed a map
$$\xymatrix{
{ K(k) } \ar[r]^{i_*} \ar[d] &
{ K(V) } \ar[r]^{j^*} \ar[d] & 
{ K(K) } \ar[r]^{\partial} \ar[d] &
{ K(k)[-1] } \ar[d] \cr
{ T(k) } \ar[r]^{i_*} &
{ T(V) } \ar[r]^{j^*} &
{ T(V|K) } \ar[r]^{\partial} &
{ T(k)[-1] } \cr
}$$
from the localization sequence in algebraic $K$-theory to an analogous
cofibration sequence relating the topological Hochschild spectra
$T(k)$ and $T(V)$ and a new topological Hochschild spectrum $T(V|K)$.
The circle group $\mathbb{T}$ acts on the terms of the lower sequence,
and one defines
$$\TR^n(V|K;p) = T(V|K)^{C_{p^{n-1}}}$$
to be the fixed points of the subgroup $C_{p^{n-1}} \subset
\mathbb{T}$ of the indicated order. As $n \geq 1$ varies, these
spectra form a pro-ring-spectrum whose structure map 
$$R \colon \TR^n(V|K;p) \to \TR^{n-1}(V|K;p)$$
is called the restriction map. We proved in \emph{op.~cit.}~that the
homotopy groups
$$\TR_q^n(V|K;p) = \pi_q(\TR^n(V|K;p))$$
of this pro-ring-spectrum form a Witt complex over the ring $V$
endowed with the canonical log-structure given by the canonical
inclusion
$$\alpha \colon M = K^* \cap V \hookrightarrow V.$$
We further showed that there exists an initial Witt complex
$\smash{W_n\Omega_{(V,M)}^q}$ over the log-ring $(V,M)$ (denoted by
$W_n\omega_{(V,M)}^q$ in \emph{op.~cit.}) and that the canonical map
$$W_n\Omega_{(V,M)}^q \to \TR_q^n(V|K;p)$$
is an isomorphism, for $q \leq 2$. The main result
of~\emph{op.~cit.}~was the evaluation, up to pro-isomorphism, of the
homotopy groups with $\Z/p^v$-coefficients
$$\TR_q^n(V|K;p,\Z/p^v) = \pi_q(\TR^n(V|K;p),\Z/p^v)$$
in the case where the field $K$ contains the $p^v$th roots of
unity. Indeed, we showed in~\emph{op.~cit.},~Thm.~C, that, as $n \geq
1$ varies, the canonical maps
$$W_n\Omega_{(V,M)}^q \otimes S_{\Z/p^v}(\mu_{p^v}) \to
\TR_q^n(V|K;p,\Z/p^v)$$
that take a generator $\zeta \in \mu_{p^v}$ to the image by the
cyclotomic trace maps of the associated Bott element $\beta_{\zeta}
\in K_2(K,\Z/p^v)$ form an isomorphism of pro-abelian groups. However,
if the field $K$ does not contain the $p^v$th roots of unity, the
structure of the pro-abelian groups $\TR_q\cs(V|K;p,\Z/p^v)$ is
presently unknown.

In this paper we consider the colimit spectra
$$\TR^n(\bar{V}|\bar{K};p) = \colim \TR^n(V_{\alpha}|K_{\alpha};p)$$
where, on the right-hand side, $K_{\alpha}$ ranges over all finite field
extension of $K$ that are contained in an algebraic closure $\bar{K}$
of $K$, and where $V_{\alpha}$ and $\bar{V}$ are the integral closures
of $V$ in $K_{\alpha}$ and $\bar{K}$, respectively. The canonical map
$$W_n\Omega_{(\bar{V},\bar{M})}^q \to \TR_q^n(\bar{V}|\bar{K};p)$$
is an isomorphism, for $q \leq 2$, but the statement for the higher
homotopy groups with $\Z/p^v$-coefficients is not valid for $\bar{K}$
as isomorphisms of pro-abelian groups are generally not preserved
under infinite colimits. The purpose of this paper is to completely
determine the $p$-adic homotopy groups
$$\TR_q^n(\bar{V}|\bar{K};p,\Zp) =
\pi_q(\TR^n(\bar{V}|\bar{K};p),\Zp).$$
We recall that, for any spectrum $X$, the $p$-adic homotopy groups are
related to the integral homotopy groups by a natural short-exact
sequence
$$0 \to \Ext(\Qp/\Zp,\pi_q(X)) \to \pi_q(X,\Zp) \to
\Hom(\Qp/\Zp,\pi_{q-1}(X)) \to 0.$$
The right-hand term is also written $T_p(\pi_{q-1}(X))$ and called the
$p$-primary Tate module of the group $\pi_{q-1}(X)$. It is not
difficult to see that $\smash{W_n\Omega_{(\bar{V},\bar{M})}^q}$ is a
divisible group, for $q$ positive, and hence we obtain a canonical
isomorphism
$$\TR_2^n(\bar{V}|\bar{K};p,\Zp) \xto{\sim}
T_p(\TR_1^n(\bar{V}|\bar{K};p)) \xleftarrow{\sim}
T_p(W_n\Omega_{(\bar{V},\bar{M})}^1).$$
This, in turn, induces a map of graded
$\TR_0^n(\bar{V}|\bar{K};p,\Zp)$-algebras
$$S_{\TR_0^n(\bar{V}|\bar{K};p,\Zp)}(T_p\TR_1(\bar{V}|\bar{K};p)) \to
\TR_*^n(\bar{V}|\bar{K};p,\Zp)$$
by using the multiplicative structure on the right-hand side.

\begin{bigthm}\label{gradedringstructure}The group
$\TR_q^n(\bar{V}|\bar{K};p)$ is divisible, if $q>0$, and uniquely
divisible, if $q>0$ and even. The $p$-primary Tate module
$T_p\TR_1^n(\bar{V}|\bar{K};p)$ is a free module of rank one over
$\TR_0^n(\bar{V}|\bar{K};p,\Zp)$, and the canonical map
$$S_{\TR_0^n(\bar{V}|\bar{K};p,\Zp)}(T_p\TR_1^n(\bar{V}|\bar{K};p))
\to \TR_*^n(\bar{V}|\bar{K};p,\Zp)$$
is an isomorphism.
\end{bigthm}

We note the formal analogy with the following result by
Suslin~\cite{suslin,suslin1} on the algebraic $K$-theory of the field
$\bar{K}$. The group $K_q(\bar{K})$ is divisible, if $q>0$, and
uniquely divisible, if $q > 0$ and even; the Tate module
$T_pK_1(\bar{K},\Zp)$ is a free module of rank one over
$K_0(\bar{K},\Zp)$, and the canonical map
$$S_{K_0(\bar{K},\Zp)}(T_pK_1(\bar{K})) \to K_*(\bar{K},\Zp)$$
is an isomorphism.

The ring $\TR_0^n(\bar{V}|\bar{K};p,\Zp)$ is canonically isomorphic to
the ring $W_n(\bar{V})\f$ given by the $p$-completion of the ring of
($p$-typical) Witt vectors in $\bar{V}$. We determine the structure of
this ring following the work of Fontaine~\cite{fontaine}. To state the
result, we first let $R_{\bar{V}}$ be the ring given by the limit of
the diagram
$$\bar{V}/p\bar{V}  \xleftarrow{\varphi} \bar{V}/p\bar{V}
\xleftarrow{\varphi} \bar{V}/p\bar{V} \xleftarrow{\varphi} \cdots$$
with the Frobenius as structure map. The ring $R_{\bar{V}}$ is a
perfect $\Fp$-algebra and an integrally closed domain whose quotient
field is algebraically closed. We show that there is a surjective ring
homomorphism
$$\theta_n\: W(R_{\bar{V}}) \twoheadrightarrow W_n(\bar{V})\f$$ 
whose kernel is a principal ideal with a generator given as
follows. We choose a sequence $\e = \{\e^{(v)}\}_{v\geq 1}$ of
compatible primitive $p^{v-1}$th roots of unity in $\bar{V} \subset
\bar{K}$. The sequence $\e$ determines, by reduction modulo $p$, an
element of $R_{\bar{V}}$ that we also denote by $\e$. We let $\e_n \in
R_{\bar{V}}$ be the unique $p^n$th root of $\e$, and let $[\e_n] \in
W(R_{\bar{V}})$ be the Teichm\"{u}ller representative. Then
$([\e]-1)/([\e_n]-1)$ generates the kernel of $\theta_n$. Moreover,
the maps $\theta_n$ satisfy that $R \circ \theta_n = \theta_{n-1}$ and
$F \circ \theta_n = \theta_{n-1} \circ F$. 

We give a similar complete description of the structure of the
$p$-primary Tate module $T_p\TR_1^n(\bar{V}|\bar{K};p)$. The
$K$-theory Bott element $\beta_{\e} = \beta_{\e}^K$ which is defined
to be the image of $\e = \{\e^{(v)}\}_{v \geq 1}$ by the canonical
isomorphism
$$T_p(\bar{K}^*) \xto{\sim}
T_pK_1(\bar{K})$$
is a $K_0(\bar{K},\Zp)$-module generator. But the $\TR$-theory Bott
element $\beta_{\e,n} = \beta_{\e,n}^{\TR}$ which is defined to be the
image of $\beta_{\e}$ by the cyclotomic trace map
$$T_pK_1(\bar{K}) \to T_p\TR_1^n(\bar{V}|\bar{K};p)$$
is not a generator of the $\TR_0^n(\bar{V}|\bar{K};p,\Zp)$-module on
the right-hand side. Instead we have the following result.

\begin{bigthm}\label{wittcomplexstructure}Let $\e =
\{\e^{(v)}\}_{v \geq 1}$ be a compatible sequence of
$p^{v-1}$th roots of unity in $\bar{V} \subset \bar{K}$. Then the
following statements hold. 

(i) The $W_n(\bar{V})\f$-module $T_p\TR_1^n(\bar{V}|\bar{K};p)$ has
a generator $\alpha_{\e,n}$ such that
$$\beta_{\e,n}=\theta_n([\e_n]-1)\cdot\alpha_{\e,n}.$$

(ii) The restriction and Frobenius maps
$$R,F \colon T_p\TR_1^n(\bar{V}|\bar{K};p) \to
T_p\TR_1^{n-1}(\bar{V}|\bar{K};p)$$
take $\alpha_{\e,n}$ to $\theta_n(([\e_{n-1}]-1)/([\e_n]-1)) \cdot
\alpha_{\e,n-1}$ and $\alpha_{\e,n-1}$, respectively.

(iii) The action of the Galois group $G_K =
\operatorname{Gal}(\bar{K}/K)$ is given by
$$\alpha_{\e,n}^{\sigma} = \chi(\sigma) \cdot
\theta_n(([\e_n]-1)/([\e_n^{\sigma}]-1)) \cdot\alpha_{\e,n},$$
where $\chi\:G_K \to \operatorname{Aut}(\mu_{p^{\infty}})=\Zp^*$ is
the cyclotomic character.
\end{bigthm}

We remark that parts~(ii) and~(iii) of Thm.~\ref{wittcomplexstructure}
are easy consequences of part~(i) and the fact that
$R(\beta_{\e,n}) = F(\beta_{\e,n}) = \beta_{\e,n-1}$. Since in
positive degrees the integral homotopy groups
$\TR_q^n(\bar{V}|\bar{K};p)$ are divisible we obtain a canonical
isomorphism
$$\TR_q^n(\bar{V}|\bar{K};p,\Qp/\Zp) \xto{\sim}
\TR_q^n(\bar{V}|\bar{K};p,\Zp) \otimes \Qp/\Zp.$$
Hence Thms.~\ref{gradedringstructure} and~\ref{wittcomplexstructure}
also determine the structure of the $G_K$-modules on the left-hand
side. The following conjecture was first formulated in~\cite{h2}.

\begin{unconj}For all positive integers $q$, the canonical map
$$\TR_q\cs(V|K;p,\Qp/\Zp) \to
\TR_q\cs(\bar{V}|\bar{K};p,\Qp/\Zp)^{G_K}$$
is an isomorphism of pro-abelian groups and the higher continuous
cohomology groups of the pro-$G_K$-module
$\TR_q\cs(\bar{V}|\bar{K};p,\Qp/\Zp)$ vanish.
\end{unconj}

We use Thms.~\ref{gradedringstructure} and~\ref{wittcomplexstructure}
and a theorem of Tate~\cite{tate} to show that the rational cohomology
groups $H_{\operatorname{cont}}^i(G_K,\TR_q\cs(\bar{V}|\bar{K};p,\Qp))$  
vanish for $i\geq 0$ and $q>0$. We hope that similar methods will make
it possible to understand the structure of the cohomology groups
$H_{\cont}^i(G_K,\TR_q\cs(\bar{V}|\bar{K};p,\Qp/\Zp))$ in question.

Finally, Thms.~\ref{gradedringstructure}
and~\ref{wittcomplexstructure} determine the structure of the $p$-adic
topological cyclic homology groups $\TC_q(\bar{V}|\bar{K};p,\Zp)$.
Indeed, we obtain the following new result.

\begin{bigthm}\label{tcthm}The cyclotomic trace induces an isomorphism
$$K_q(\bar K,\Zp)\xto{\sim}\TC_q(\bar V|\bar K;p,\Zp),$$
for all integers $q$.
\end{bigthm}

We remark that the common group in the statement of Thm.~\ref{tcthm}
is canonically isomorphic to the $G_K$-module $\Zp(q/2)$, if $q$ is a 
non-negative even integer, and is zero, otherwise. We also note that
the vanishing of $\TC_q(\bar{V}|\bar{K};p,\Zp)$, for $q$ odd
(including $q = -1$), is an immediate consequence of the fact that the
quotient field of the $\Fp$-algebra $R_{\bar{V}}$ is algebraically
closed.

The paper is organized as follows. Sect.~\ref{wittvectors} determines
the structure of the ring $W_n(\bar{V})\f$; Sect.~\ref{algclosed}
contains the proofs of Thms.~\ref{gradedringstructure}
and~\ref{wittcomplexstructure} with the exception of the existence of
the generator $\alpha_{\e,n}$ which is proved in
Sect.~\ref{generatorsection}; and the final Sect.~\ref{galoissection}
concerns Galois cohomology.

The results of this paper were reported in expository form
in~\cite{h2}.

It is a great pleasure to thank the University of Tokyo and Takeshi
Saito in particular for their hospitality during the writing of
parts of this paper. The author also expresses his sincere gratitude
to an anonymous referee who read an earlier version of this paper with
extreme care and made numerous helpful suggestions.

\section{Witt vectors}\label{wittvectors}

\subsection{}In this section we determine the structure of the ring
$$W_n(\bar{V})\f \xto{\sim} \TR_0^n(\bar{V}|\bar{K};p,\Zp)$$
that is given by the $p$-completion of the ring of ($p$-typical) Witt
vectors of length $n$ in $\bar{V}$. This follows the work of
Fontaine~\cite{fontaine}. We refer the reader to~\cite[Sect.~1.1]{hm3}
for generalities on Witt vectors.

We say that a ring $A$ is $p$-complete if
the canonical map from $A$ to the $p$-completion $A\f = \lim_v A/p^vA$
is an isomorphism, or equivalently, if the $p$-adic topology on $A$ is
complete and separated.

\begin{lemma}Let $A$ be a $p$-complete and $p$-torsion free ring. Then
$W(A)$ is again a $p$-complete and $p$-torsion free ring.
\end{lemma}

\begin{proof}Since $A$ is $p$-torsion free then so is $W_n(A)$, and
hence, the sequences
$$0 \to A/p^vA \xto{V^{n-1}} W_n(A)/p^vW_n(A) \xto{R}
W_{n-1}(A)/p^vW_{n-1}(A) \to 0$$
are exact. It follows, by induction, that $W_n(A)$ is a $p$-complete
ring. Taking limits over $n$, the exact sequences
$$0\to W_n(A) \xto{p^v}W_n(A) \to W_n(A)/p^vW_n(A) \to 0$$
give rise to an isomorphism
$$W(A)/p^vW(A)\xto{\sim}\lim_n(W_n(A)/p^vW_n(A)).$$
Taking the limit over $v$ and using that limits commute, we find that
the lower horizontal map in following diagram is an isomorphism.
$$\xymatrix{
{W(A)} \ar[r] \ar[d] &
{\lim_nW_n(A)} \ar[d] \cr
{W(A)\f} \ar[r] &
{\lim_n(W_n(A)\f).} \cr
}$$
The right-hand vertical map is an isomorphism, since $W_n(A)$ is a
$p$-complete ring, and the top horizontal map is an isomorphism for
trivial reasons. Hence, the left-hand vertical map is an isomorphism
as stated.
\end{proof}

\begin{lemma}Let $A$ be a ring without $p$-torsion. Then the canonical
maps $W(A)\f \to W(A\f)\f$ and $W_n(A)\f \to W_n(A\f)\f$ are
isomorphisms.
\end{lemma}

\begin{proof}Since $A$ is $p$-torsion free, the canonical map
$A/p^vA \to A\f/p^vA\f$ is an isomorphism. Hence, inductively, so
are the induced maps
$$W_n(A)/p^vW_n(A)\xto{\sim}W_n(A\f)/p^vW_n(A\f),$$
for all $n\geq 1$. Using that $W_n(A)$ and $W_n(A\f)$ are $p$-torsion
free, we conclude that
$$W(A)/p^vW(A)\xto{\sim}W(A\f)/p^vW(A\f).$$
Hence, the induced map of limits over $v$ is an isomorphism, and this
is the statement of the lemma.
\end{proof}

\begin{lemma}\label{reductionisomorphism}Let $A$ be a $p$-complete and
$p$-torsion free ring. Then the canonical projection
$$\lim_FW(A) \to \lim_FW(A/pA)$$
is an isomorphism, and the inverse is given as follows. Let
$x=\{x^{(m)}\}_{m \geq 1}$ be an element of $\lim_F W(A/pA)$, and let
$\tilde x^{(m)} \in W(A)$ be liftings of the individual
$x^{(m)}$. Then the sequence $\{ F^k(\tilde x^{(m+k)}) \}_{k \geq 0}$
converges in the $p$-adic topology to an element $\hat{x}^{(m)}\in
W(A)$, and the sequence $\hat{x}=\{ \hat{x}^{(m)} \}_{m \geq 1}$ is
the element of $\lim_F W(A)$ whose image by the canonical projection
is $x$.
\end{lemma}

\begin{proof}We first show that the canonical projection is an
isomorphism. To this end, we consider the six-term exact sequence
associated with the following short-exact sequence of limit systems.
$$0\to W(pA)\to W(A)\to W(A/pA)\to 0.$$
It follows that it suffices to show that $\lim_FW(pA)$ and
$R^1\lim_F W(pA)$ vanish. These groups are given by the kernel and
cokernel, respectively, of the map
$$\id-\sigma\colon \prod_{m \geq 1} W(pA) \to \prod_{m \geq 1} W(pA),$$
where $\sigma$ maps the $m$th factor to the $(m-1)$st factor by the
Frobenius. Since $A$ is $p$-torsion free, the ghost map defines an
isomorphism
$$w\colon W(pA) \xto{\sim} \prod_{s\geq 0} p^{s+1}A.$$
So the map in question is isomorphic to the map
$$\id-\sigma \colon \prod p^{s+1}A \to \prod p^{s+1}A$$
where the products ranges over $s \geq 0$ and $m \geq 1$ and where
$\sigma$ maps the $(s,m)$th factor to the $(s-1,m-1)$st factor
by the canonical inclusion $p^{s+1}A\to p^sA$. Hence, $\id-\sigma$ is
an isomorphism with inverse
$$(\id-\sigma)^{-1}=\sum_{n\geq 0}\sigma^n.$$
The series on the right converges since $A$ is $p$-complete.

It remains to show that the inverse of the canonical projection is
given as stated. We must show that the sequence
$\{ F^k(\tilde{x}^{(m+k)}) \}_{k \geq 0}$ converges with a unique limit,
and since the $p$-adic topology on $W(A)$ is complete and 
separated, it suffices to show that the sequence is Cauchy. To this
end, we first note that a sequence $\{x_s\}_{s \geq 0}$ is in the
image of the ideal $W(pA)$ by the ghost map
$$w\colon W(A)\to A^{\N_0}$$
if and only if $x_s\in p^{s+1}A$, for all $s \geq 0$. Now
$$\aligned
w_s(F^k(\tilde x^{(m+k)}))&=w_{s+k}(\tilde x^{(m+k)})\cr
{}&=(\tilde x_0^{(m+k)})^{p^{s+k}}+p(\tilde x_1^{(m+k)})^{p^{s+k-1}}+\dots
+p^{s+k}\tilde x_{s+k}^{(m+k)}\cr
\endaligned$$
which shows that
$$w_s(F^{k+1}(\tilde x^{(m+k+1)}))-w_s(F^k(\tilde x^{(m+k)}))\in
p^{s+k+1}A.$$
Hence the sequence $\{F^k(\tilde x^{(m+k)})\}_{k \geq 0}$ is Cauchy as
desired.
\end{proof}

\begin{lemma}\label{finitewitt}Let $A$ be a $p$-complete and
$p$-torsion free ring. Then the canonical projection induces an
isomorphism
$$\lim_FW(A) \xto{\sim} \lim_FW_n(A).$$
\end{lemma}

\begin{proof}We consider the exact sequence of limit systems with
$n$th terms
$$0 \to W(A) \xto{V^n} W(A) \xto{\pr} W_n(A) \to 0.$$
The structure maps in the middle and right-hand limit systems are
given by the Frobenius, and the structure map in the left-hand limit
system is given by multiplication by $p$. The induced six-term exact
sequence takes the form
$$\begin{aligned}
0 & \to \lim_p W(A) \to \lim_F W(A) \to \lim_F W_n(A) \cr
{} & \to R^1\lim_p W(A) \to R^1\lim_F W(A) \to R^1\lim_F W_n(A) \to 0.
\end{aligned}$$
The first and fourth terms vanish since the $p$-adic topology on
$W(A)$ is separated and complete, respectively. The lemma follows.
\end{proof}

\subsection{}Let $V$ be a complete discrete valuation ring with
quotient field $K$ and perfect residue field $k$ of mixed
characteristic $(0,p)$. Let $\bar{K}$ be an algebraic closure of $K$
and let $\bar{V}$ 
be the integral closure of $V$ in $\bar{K}$. Let
$$\bar{V}\f=\lim_n\bar{V}/p^n\bar{V}$$
be the $p$-completion of $\bar{V}$, and let $\bar{K}\f$ be the
quotient field of $\bar{V}\f$. Then $\bar{K}\f$ is again algebraically
closed~\cite[Chap.~2,~Thm.~12]{artin}. The valuation on $K$ extends
uniquely to a valuation on $\bar{K}$ (resp.~$\bar{K}\f$) with value
group the additive group of rational numbers and $\bar{V}$
(resp.~$\bar{V}\f$) is the valuation ring. In particular, $\bar{V}$
(resp.~$\bar{V}\f$) is an integrally closed local ring of dimension
one. However, $\bar{V}$ (resp.~$\bar{V}\f$) is not a
noetherian ring. We normalize the valuation on $\bar{V}$ (resp.~$\bar
V\f$) such that $v(p)=1$. 

Let $x,y\in\bar{V}$ and suppose that $x^p=y^p$. Then
$v(x-y)\geq 1/(p-1)$. Indeed, we have $x=\zeta y$ with $\zeta^p=1$,
and hence, $x-y=(\zeta-1)y$. But $\zeta-1$ is either zero or a
uniformizer of $\Qp(\mu_p)$, and thus, $v(\zeta-1)\geq 1/(p-1)$.

We follow~\cite[Sect.~1.1]{hm3} and write $[a]_n \in W_n(A)$ for
the Teichm\"{u}ller representative of $a \in A$. We note that the
notation $\underline{a}_n$ was used instead of $[a]_n$ in~\cite{hm2}.

\begin{lemma}\label{modpwitt}The map $\theta_n' \colon \bar{V}\to
W_n(\bar V)/pW_n(\bar{V})$ given by $\theta_n'(x)=[x]_n^p$ is a
surjective ring homomorphism whose kernel is the ideal of elements of
valuation greater than or equal to $(1-p^{-n})/(p-1)$. Moreover,
$R(\theta_n'(x))=\theta_{n-1}'(x)$,
$\smash{F(\theta_n'(x))=\theta_{n-1}'(x^p)}$, and
$\smash{V(\theta_{n-1}'(x))=\theta_n'((-px)^{p^{-1}})}$.
\end{lemma}

\begin{proof}We recall from \cite[Lemma~3.1.2]{hm2} that for every
ring $A$, and for all $x,y\in A$, we have
$$[x]_n^p+[y]_n^p=([x]_n+[y]_n)^p=([x+y]_n)^p$$
in $W_n(A)/pW_n(A)$. Moreover, if $p$ is odd, then in addition
$$V(1)=[-p]_n$$ in $W_n(A)/pW_n(A)$ by
\emph{op.~cit.},~Lemma~3.1.1. It follows that the map $\theta_n'$ is a
ring homomorphism and that $\theta_n'((-p)^{p^{-1}})=V(1)$. By easy
induction, we find
$$\theta_n'(a^{p^{-(i+1)}}(-p)^{(1-p^{-i})/(p-1)})
=[a]_n^{p^{-i}}V^i(1)=V^i([a]_{n-i}),$$
which shows that $\theta_n'$ is surjective. Moreover,
$\theta_n'(x)\in V^iW_n(\bar{V})/pV^iW_n(\bar{V})$ if and only if
$v(x)\geq(1-p^{-i})/(p-1)$.
\end{proof}

\begin{cor}\label{frobsurj}The Frobenius $F\colon W_n(\bar{V})\f\to
W_{n-1}(\bar{V})\f$ is surjective.
\end{cor}

\begin{proof}It follows from Lemma~\ref{modpwitt} that the statement
holds after reduction modulo $p$. Let $W_{n,v}(\bar{V}) = 
W_n(\bar{V})/p^vW_n(\bar{V})$. Then an induction argument based on the
diagram
$$\xymatrix{
{0} \ar[r] &
{W_{n,1}(\bar{V})} \ar[r]^{p^{v-1}} \ar[d]^{F} &
{W_{n,v}(\bar{V})} \ar[r] \ar[d]^{F} &
{W_{n,v-1}(\bar{V})} \ar[r] \ar[d]^{F} &
{0} \cr
{0} \ar[r] &
{W_{n-1,1}(\bar{V})} \ar[r]^{p^{v-1}} &
{W_{n-1,v}(\bar{V})} \ar[r] &
{W_{n-1,v-1}(\bar{V})} \ar[r] &
{0} \cr
}$$
shows that the statement is true after reduction modulo $p^v$, for all
$v\geq 1$. Since the vertical maps in the diagram above are all
surjective, the sequence of kernels
$$0\to K_{n,1}\to K_{n,v}\to K_{n,v-1}\to 0$$
is exact. Hence, the sequence of limit systems
$$0\to K_{n,v}\to W_{n,v}(\bar{V})\xto{F}W_{n-1,v}(\bar{V})\to 0$$
gives an exact sequence
$$
0\to\lim_vK_{n,v}
\to\lim_vW_{n,v}(\bar{V})
\xto{F}\lim_vW_{n-1,v}(\bar{V})\to 0.$$
The derived limit $R^1\lim_v K_{n,v}$ vanishes since the structure maps
are surjective.
\end{proof}

Following Fontaine \cite{fontaine} we define $R_{\bar V}$ to be the
limit of the diagram
$$\dots\to\bar{V}/p\bar{V}\xto{\varphi}\bar{V}/p\bar{V}
\xto{\varphi}\bar{V}/p\bar{V},$$
where $\varphi$ is the Frobenius. We note that the map of limits
induced from the following map of towers of multiplicative monoids is
an isomorphism.
$$\xymatrix{
{\dots} \ar[r] &
{\bar{V}\f} \ar[r] \ar[d] &
{\bar{V}\f} \ar[r] \ar[d] &
{\bar{V}\f} \ar[d] \cr
{\dots} \ar[r] &
{\bar{V}/p\bar{V}} \ar[r]^{\varphi} &
{\bar{V}/p\bar{V}} \ar[r]^{\varphi} &
{\bar{V}/p\bar{V}.} \cr
}$$
The horizontal maps are given by raising to the $p\,$th
power. Indeed, given an element $x=\{x^{(v)}\}_{v\geq 1}\in
R_{\bar{V}}$ we choose liftings  $\tilde x^{(v)} \in \bar{V}$ of
$x^{(v)}$. Then the sequence $\{(\tilde x^{(m+v)})^{p^m}\}_{m\geq 0}$
converges in $\bar{V}\f$ and the limit
$$\hat x^{(v)}=\lim_{m\to\infty}(\tilde x^{(m+v)})^{p^m}$$
is independent of the choice of liftings $\tilde x^{(v)}$ of
$x^{(v)}$. We define a valuation $v_R$ of $R_{\bar{V}}$ by
$v_R(x)=v_{\bar{V}\f}(\hat x^{(1)})$. We note that the kernel of the
canonical projection
$$\pr_1 \colon R_{\bar{V}} \to \bar{V}/p\bar{V}$$
is equal to the ideal of elements of valuation greater than or equal
to $1$.



\begin{prop}\label{thewittring}There is a natural surjective ring
homomorphism
$$\theta_n \colon W(R_{\bar{V}}) \twoheadrightarrow W_n(\bar{V})\f$$
whose kernel is a principal ideal. If $\e=\{\e^{(v)}\}_{v\geq 1}$ is a
compatible sequence of primitive $p^{v-1}$th roots of unity in
$\bar{V}$ considered
as an element of $R_{\bar{V}}$, and if $\e_n$ is the unique $p^n$th
root of $\e$, then $([\e]-1)/([\e_n]-1)$ generates the kernel of
$\theta_n$. Moreover, the maps $\theta_n$ satisfy that $R \circ
\theta_n = \theta_{n-1}$ and $F \circ \theta_n = \theta_{n-1} \circ F$.
\end{prop}

\begin{proof}We consider the ring isomorphism
$$\psi \colon W(R_{\bar{V}}) \to \lim_F W_n(\bar{V})\f$$
defined to be the composite of the following isomorphisms.
$$W(R_{\bar{V}}) \xto{\sim}
\lim_FW(\bar{V}/p\bar{V}) \xleftarrow{\sim}
\lim_FW(\bar{V})\f \xto{\sim}
\lim_FW_n(\bar{V})\f.$$
The first map is an isomorphism since the underlying set of the Witt
ring $W(A)$ is the product of copies of $A$ indexed by the set of
non-negative integers; the second map is an isomorphism by
Lemma~\ref{reductionisomorphism}; and the third map is an isomorphism
by Lemma~\ref{finitewitt}. The map of the statement is then induced
from the composite
$$W(R_{\bar{V}}) \xto{F^{n-1}} 
W(R_{\bar{V}}) \xto{\psi}
\lim_FW_n(\bar{V})\f \xto{\pr_n}
W_n(\bar{V})\f.$$
We first show that $\theta_n(\xi_n) = 0$. The map $F^{n-1}$ takes
$[\e_m]$ to $[\e_{m-(n-1)}]$, and the canonical isomorphism
$$W(R_{\bar{V}}) = W(\lim_{\varphi} \bar{V}/p\bar{V}) \xto{\sim}
\lim_F W(\bar{V}/p\bar{V})$$
takes $[\e_k]$ to the sequence $\{[\e_k^{(v)}]\}_{v \geq 1}$. To find
the image of this sequence by the inverse of the canonical isomorphism
$$\lim_F W(\bar{V}/p\bar{V}) \xleftarrow{\sim} \lim_F W(\bar{V})\f,$$
we use the formula given by Lemma~\ref{reductionisomorphism} but with
a particular choice of liftings $\tilde{x}^{(v)}$ to $W(\bar{V})\f$ of
the elements $[\e_k^{(v)}]$ of $W(\bar{V}/p\bar{V})$. Let
$\tilde{\e}_k = \{\tilde{\e}_k^{(v)}\}_{v \geq 1}$ be a sequence of
elements in $\bar{V}$ such that $\e_k^{(v)}$ is a primitive
$p^{k+v-1}$th root of unity, for $k + v - 1$ non-negative, and such
that $(\e_k^{(v+1)})^p = \e_k^{(v)}$, for all $v \geq 1$. Then we
choose $\tilde{x}_k^{(v)} = [\tilde{\e}_k^{(v)}]$. Since these
elements are already compatible under the Frobenius map, we find that
also $\hat{x}_k^{(v)} = [\tilde{\e}_k^{(v)}]$. Hence, the image of the
sequence $\{[\e_k^{(v)}]\}_{v \geq 1}$ by the inverse of the
isomorphism above, is the sequence $\{[\tilde{\e}_k^{(v)}]\}_{v \geq 1}$. 
But this sequence is mapped to $1$ by the composite
$$\lim_F W(\bar{V})\f \xto{\sim} \lim_F W_n(\bar{V})\f \xto{\pr}
W_n(\bar{V})\f$$
if and only if $k \leq -(n-1)$. In particular, $\theta_n(\xi_n)$ is
zero as stated.

It remains to show that $\xi_n$ generates the kernel of $\theta_n$. We
show inductively that for all $v\geq 1$, the induced map
$$\theta_{n,v} \colon W_v(R_{\bar{V}})/(\xi_{n,v}) \to
W_n(\bar{V})/p^vW_n(\bar{V})$$
is an isomorphism. Here $\xi_{n,v} = ([\e]_v - 1)/([\e_n]_v - 1)$ is
the image by the canonical projection of $\xi_n$ in
$W_v(R_{\bar{V}})$. In the basic case $v = 1$, we have the following
commutative diagram.
$$\xymatrix{
{R_{\bar{V}}} \ar[r]^(.38){\pr} \ar[d]^{\varphi^{-1}} &
{R_{\bar{V}}/(\xi_{n,1})} \ar[r]^(.4){\theta_{n,1}} &  
{W_n(\bar{V})/pW_n(\bar{V})} \ar@{=}[d] \cr
{R_{\bar{V}}} \ar[r]^(.45){\pr_1} &
{\bar{V}/p\bar{V}} \ar[r]^(.35){\bar{\theta}_n'} &
{W_n(\bar{V})/pW_n(\bar{V}).} \cr
}$$
By Lemma~\ref{modpwitt}, the kernel of the composition of the maps in
the lower row is equal to the ideal of elements of valuation greater
than or equal to $(1-p^{-n})/(p-1)$. The element $\xi_{n,1}$, which
generates the kernel of the composition of the maps in the top two,
has valuation
$$v_R(\xi_{n,1}) = v_R(\e - 1) - v_R(\e_n - 1) = \frac{1}{p^{-1}(p -
  1)} - \frac{1}{p^{n-1}(p-1)} = \frac{1 - p^{-n}}{1 - p^{-1}}.$$
Since $\varphi^{-1}$ is an isomorphism and $v_R(\varphi^{-1}(x)) =
p^{-1} v_R(x)$, it follows that $\theta_{n,1}$ is an isomorphism as
desired.

Finally, to prove the induction step it suffices to show that the top
row in the following diagram is exact.
$$\xymatrix{
{0} \ar[r] &
{R_{\bar{V}}/(\xi_{n,1})} \ar[r]^(.45){p^{v-1}}
\ar[d]^{\theta_{n,1}}_{\sim} & 
{W_v(R_{\bar{V}})/(\xi_{n,v})} \ar[r]^(.45){R} \ar[d]^{\theta_{n,v}} &
{W_{v-1}(R_{\bar{V}})/(\xi_{n,v-1})} \ar[r] \ar[d]^{\theta_{n,v-1}}_{\sim} &
{0} \cr
{0} \ar[r] &
{W_n(\bar{V})/p} \ar[r]^(.47){p^{v-1}} &
{W_n(\bar{V})/p^v} \ar[r] &
{W_n(\bar{V})/p^{v-1}} \ar[r] &
{0} \cr
}$$
But the sequence
$$\xymatrix{
{0} \ar[r] &
{R_{\bar{V}}} \ar[r]^(.4){p^{n-1}} &
{W_v(R_{\bar{V}})} \ar[r]^(.47){R} &
{W_{v-1}(R_{\bar{V}})} \ar[r] &
{0} \cr
}$$
is exact, since $R_{\bar{V}}$ is a perfect $\Fp$-algebra, and the
element $\xi_{n,1} \in R_{\bar{V}}$ is a non-zero-divisor. It follows,
inductively, that $\xi_{n,v} \in W_v(R_{\bar{V}})$ is a non-zero-divisor
and that the sequence in question is exact.
\end{proof}

\begin{addendum}\label{thelimitwittring}There is a ring isomorphism
$$\psi \colon W(R_{\bar{V}}) \xto{\sim} \lim_FW_n(\bar{V})\f$$
given by $\psi(a) = \{\theta_n(F^{-(n-1)}(a))\}_{n\geq 1}$. Moreover,
$R \circ \psi = \psi \circ F^{-1}$.
\end{addendum}

\begin{proof}We showed at the beginning of the proof of
  Prop.~\ref{thewittring} that the map $\psi$ is an isomorphism. The
  statement that $R \circ \psi = \psi \circ F^{-1}$ follows easily
  from the formula $R \circ \theta_n = \theta_{n-1}$.
\end{proof}


\subsection{}We conclude this section with some miscellaneous results
about the Witt rings $W(R_{\bar{V}})$ and $W_n(\bar{V})\f$. The
following result is~\cite[1.3.1]{fontaine}. 

\begin{prop}\label{rv}The ring $R_{\bar{V}}$ is a perfect
\,$\Fp$-algebra and an integrally closed domain whose quotient field
is algebraically closed and complete with respect to the valuation
$v_R$ with valuation ring $R_{\bar{V}}$. \hfill\space\qed
\end{prop}


\begin{addendum}\label{domain}The ring $W(R_{\bar{V}})$ is an integral
domain. 
\end{addendum}

\begin{proof}It suffices, by Prop.~\ref{rv}, to show that if $A$ is an
$\Fp$-algebra and an integral domain, then $W(A)$ is an integral
domain. Suppose, conversely, that $W(A)$ is not an integral domain. We
show that we can find a zero-divisor $x\in W(A)$ such that $x+VW(A)$
is a zero-divisor in $A=W(A)/VW(A)$. Let $x,y\in W(A)$ be non-zero
and suppose $xy=0$. We write $x=V^s(x')$ and $y=V^t(y')$ with $s$ and
$t$ as large as possible. Assuming that $s\geq t$,
$$xy=V^s(x')V^t(y')=V^s(x'F^sV^t(y'))=p^tV^s(x'F^{s-t}(y')).$$
We claim that $W(A)$ has no $p$-torsion. Granting this, we see that
$x'+VW(A)$ is a zero-divisor in $A=W(A)/VW(A)$. It remains to prove 
the claim. Since $A$ is an $\Fp$-algebra, we have $p=V(1)$, and hence,
for $a\in W(A)$, $pa=V(1)a=V(F(a))$. We must show that $F\colon W(A)\to
W(A)$ is a monomorphism. Again since $A$ is an $\Fp$-algebra, $F$ is
induced from the Frobenius $\varphi\colon A\to A$. But $\varphi \colon
A \to A$ is a monomorphism since $A$ is a domain.
\end{proof}

\begin{cor}\label{calctc}Let $\e=\{\e^{(v)}\}_{v\geq 1}$ be a
compatible sequence of primitive $p^{v-1}$th roots of unity considered
as an element of $R_{\bar{V}}$, and let $\e_n$ be the unique $p^n$th
root of $\e$. Then for every non-negative integer $q$, the following
sequence is exact.
$$0\to W(\Fp)\xto{([\e_1]-1)^q}
W(R_{\bar{V}})\xto{1-(\frac{[\e_1]-1}{[\e_2]-1})^qW(\varphi^{-1})}
W(R_{\bar{V}})\to 0.$$
\end{cor}

\begin{proof}Let $a =([\e_1]_n - 1)/([\e_2]_n - 1)$. Then the zeroth
Witt coordinate is $a_0 = a(0) = (\e_1 - 1)/(\e_2 - 1)$. We show by
induction on $n\geq 1$ that the sequences
$$\xymatrix{
{0} \ar[r] &
{W_n(\Fp)} \ar[rr]^{[\e_1]_n-1} &&
{W_n(R_{\bar{V}})} \ar[rr]^{1-aW_n(\varphi^{-1})} &&
{W_n(R_{\bar{V}})} \ar[r] &
{0} \cr
}$$
are exact. The basic case $n=1$ and the induction step are
similar. Indeed, it will suffice to show that for all $n\geq 1$, the
top horizontal sequence in the following diagram is exact.
$$\xymatrix{
{} &
{0} \ar[d] &&
{0} \ar[d] &&
{0} \ar[d] &
{} \cr
{0} \ar[r] &
{\Fp} \ar[rr]^{(\e_1)^{p^{n-1}}-1} \ar[d]^{V^{n-1}} &&
{R_{\bar{V}}} \ar[rr]^{1-a_0^{p^{n-1}}\varphi^{-1}} \ar[d]^{V^{n-1}} &&
{R_{\bar{V}}} \ar[d]^{V^{n-1}} \ar[r] &
{0} \cr
{0} \ar[r] &
{W_n(\Fp)} \ar[rr]^{[\e_1]_n-1} \ar[d]^{R} &&
{W_n(R_{\bar{V}})} \ar[rr]^{1-aW_n(\varphi^{-1})} \ar[d]^{R} &&
{W_n(R_{\bar{V}})} \ar[r] \ar[d]^{R} &
{0} \cr
{0} \ar[r] &
{W_{n-1}(\Fp)} \ar[rr]^{[\e_1]_{n-1}-1} \ar[d] &&
{W_{n-1}(R_{\bar{V}})} \ar[rr]^{1-aW_{n-1}(\varphi^{-1})} \ar[d] &&
{W_{n-1}(R_{\bar{V}})} \ar[r] \ar[d] &
{0} \cr
{} &
{0} &&
{0} &&
{0} &
{} \cr
}$$
We mention that the upper right-hand square commutes since
$$aV^{n-1}(x) = V^{n-1}(F^{n-1}(a)x) = V^{n-1}(a_0^{p^{n-1}}x)$$
and similar for the upper left-hand square. To prove that the top
horizontal sequence in the diagram above is exact we compose the
right hand map in the sequence with the isomorphism $\varphi\colon
R_{\bar{V}}\to R_{\bar{V}}$. The composite map takes $x$ to
$x^p-a_0^{p^n}x$. Since $R_{\bar{V}}$ is an integrally closed domain
whose quotient field is algebraically closed, the exactness follows.
\end{proof}

\begin{lemma}\label{rootsofunity}The $p^m$th roots of unity in the
ring\, $W_n(\bar{V})\f$ are the $p^m$ elements of the form $[\zeta]_n$
where $\zeta$ is a $p^m$th root of unity in the ring $\bar{V}\f$.
\end{lemma}

\begin{proof}The number of $p^m$th roots of unity in the ring
  $\bar{V}\f$ is equal to $p^m$. Hence it suffices to show that every
$p^m$th root of unity in $W_n(\bar{V})\f$ is of the form $[\zeta]_n$
where $\zeta$ is a $p^m$th root of unity in $\bar{V}\f$. Now an
element $a \in W_n(\bar{V})\f$ is a $p^m$th root of unity if an only
if the image of $a$ by the ghost map
$$w \colon W_n(\bar{V})\f \to (\bar{V}\f)^n$$
is a $p^m$th root of unity. Indeed, the ghost map is an injective ring
homomorphism since the ring $\bar{V}\f$ has no $p$-torsion. This, in
turn, is equivalent to the statement that for all $0 \leq i < n$, the
ghost coordinates
$$w_i = a_0^{p^i} + pa_1^{p^{i-1}} + \dots + p^{i-1}a_{i-1}^p +
p^ia_i$$
are $p^m$th roots of unity in $\bar{V}\f$. We show by induction on $0
\leq i < n$ that $a_0$ is a $p^m$th root of unity in $\bar{V}\f$ and
that $a_1,\dots,a_i$ are equal to zero. The statement for $i = 0$ is
trivial since $a_0 = w_0$. So we assume that $a_0$ is a $p^m$th root
of unity and that $a_1 = \dots = a_{i-1} = 0$ and must show that also
$a_i = 0$. The defining equation for the $i$th ghost coordinate yields
the equation
$$w_i - a_0^{p^i} = p^ia_i.$$
Since $w_i$ and $a_0$ are both $p^m$th roots of unity, the left-hand
side is either zero or an element of valuation at most
$1/(p-1)$. Similarly, the right-hand side is either zero or an element
of valuation at least $i$. It follows that $a_i = 0$ as desired.
\end{proof}

\begin{remark}In the proof of Lemma~\ref{rootsofunity} we have used
that $p$ is an odd prime. If $p = 2$, the $2^m$th roots of unity in
$W_n(\bar{V})\f$ are instead the elements of the form $\pm [\zeta]_n$
where $\zeta$ is a $2^m$th root of unity in $\bar{V}\f$. Hence, the
number of $2^m$th roots of unity in $W_n(\bar{V})\f$ is equal to
$2^m$, if $n = 1$, and $2^{m+1}$, if $n > 1$.
\end{remark}

\section{The algebraically closed case}\label{algclosed}

\subsection{}In this section we determine the structure of the groups
$$\TR_q^n(\bar{V}|\bar{K};p,\Zp) =
\pi_q(T(\bar{V}|\bar{K})^{C_{p^{n-1}}},\Zp)$$
and prove Thms.~\ref{gradedringstructure}
and~\ref{wittcomplexstructure} of the introduction. We first extend
some of the results from~\cite{hm2} to the algebraically closed case.

The $p$-completion and $p$-cocompletion of a spectrum $X$ are defined
by
$$\begin{aligned}
X\f & = F(M(\Qp/\Zp,-1),X) \cr
X\cf & = M(\Qp/\Zp,-1)\wedge X \cr
\end{aligned}$$
respectively~\cite{bousfield}. Here $M(\Qp/\Zp,-1)$ is a Moore
spectrum for the group $\Qp/\Zp$ concentrated in degree $-1$. We write
$\pi_q(X,\Zp)$ and $\pi_q(X,\Qp/\Zp)$ for the $q$th homotopy group of
$X\f$ and $X\cf$, respectively. There are natural short-exact
sequences
$$\begin{aligned}
{} & 0 \to \Ext(\Qp/\Zp,\pi_sX) \to \pi_s(X,\Zp) \to
\Hom(\Qp/\Zp,\pi_{s-1}X) \to 0 \cr
{} & 0 \to \pi_{s+1} X \otimes \Qp/\Zp \to \pi_s(X,\Qp/\Zp) \to
\Tor(\pi_sX,\Qp/\Zp) \to 0\cr
\end{aligned}$$
which relate the integral homotopy groups and the homotopy groups with
$\Zp$-coefficients and with $\Qp/\Zp$-coefficients. A map of spectra
induces a weak equivalence after $p$-completion if and only if it
induces an isomorphism of homotopy groups with $\Z/p$-coefficients. 

Let $\mathbb{T}$ be the circle group and let $C_r$ be the subgroup of
order $r$. We begin with the following lemma from equivariant homotopy
theory. We refer to~\cite[Sect.~4]{hm2} for the definition and properties
of the Tate spectrum.

\begin{lemma}\label{Ttate}Let $X$ be a $\mathbb{T}$-spectrum and
suppose the $\pi_tX$ vanishes for $t \ll 0$. Then the canonical map of
Tate spectra
$$\tate(\mathbb{T},X) \to \holim_n\tate(C_{p^n},X)$$
becomes a weak equivalence after $p$-completion.
\end{lemma}

\begin{proof}We show first that the map
$$(F(E_+,X)^{\mathbb{T}})\f\to\holim_n(F(E_+,X)^{C_{p^n}})\f$$
is a weak equivalence. This can be rewritten as
$$F((E_+)\cf,X)^{\mathbb{T}}\to
F((\hocolim_n\mathbb{T}/C_{p^n+}\wedge E_+)\cf,X)^{\mathbb{T}}$$
with the map induced from the canonical map
$$(\hocolim_n\mathbb{T}/C_{p^n+}\wedge E_+)\cf\to (E_+)\cf.$$
This is a weak equivalence of pointed $\mathbb{T}$-spaces if and only
if the map
$$(\hocolim_n\mathbb{T}/C_{p^n+})\cf\to (S^0)\cf$$
is a weak equivalence of pointed spaces. The latter is easily verified
by calculating homotopy groups.

We next show by induction that for $s\geq 0$, the canonical map
$$((\tilde E_s\wedge F(E_+,X))^{\mathbb{T}})\f\to
(\holim_n(\tilde E_s\wedge F(E_+,X))^{C_{p^n}})\f$$
is a weak equivalence. The basic case $s=0$ was proved above. To prove
the induction step we consider the cofibration sequence
$$\tilde E_{s-1}\rightarrowtail\tilde E_s\to\Sigma^{2s}\mathbb{T}_+
\to \Sigma\tilde{E}_{s-1}.$$
It shows that it suffices to show that the map
$$((\Sigma^{2s}\mathbb{T}_+\wedge F(E_+,X))^{\mathbb{T}})\f\to
(\holim_n(\Sigma^{2s}\mathbb{T}_+\wedge F(E_+,X))^{C_{p^n}})\f$$
is a weak equivalence. Since the inclusion $X\to F(E_+,X)$ is a
non-equivariant equivalence, we may instead show that the map
$$((\Sigma^{2s}\mathbb{T}_+\wedge X)^{\mathbb{T}})\f\to
(\holim_n(\Sigma^{2s}\mathbb{T}_+\wedge X)^{C_{p^n}})\f$$
is a weak equivalence. We now use the weak equivalence of
$\mathbb{T}$-spectra
$$\Sigma F(\mathbb{T}_+,X)\xto{\sim}\mathbb{T}_+\wedge X$$
to rewrite the latter map as
$$(\Sigma^{2s+1}F(\mathbb{T}_+,X)^{\mathbb{T}})\f\to
(\holim_n(\Sigma^{2s+1}F(\mathbb{T}_+,X)^{C_{p^n}})\f.$$
As in the basic case $s = 0$, we can rewrite this map as the map
$$\Sigma^{2s+1} F((\mathbb{T}_+)\cf,X)^{\mathbb{T}} \to
\Sigma^{2s+1} F((\hocolim_n\mathbb{T}/C_{p^n+}\wedge\mathbb{T}_+)\cf, 
X)^{\mathbb{T}}$$
induced from the canonical map
$$(\hocolim_n\mathbb{T}/C_{p^n+}\wedge\mathbb{T}_+)\cf
\to(\mathbb{T}_+)\cf.$$
We see that this is a weak equivalence of pointed $\mathbb{T}$-spaces
by calculating homotopy groups as before.

Finally, we consider the diagram
$$\xymatrix{
{((\tilde E_s\wedge F(E_+,X))^{\mathbb{T}})\f} \ar[r]^(.43){\sim} \ar[d] &
{\holim_n((\tilde E_s\wedge F(E_+,X))^{C_{p^n}})\f} \ar[d] \cr
{((\tilde E\wedge F(E_+,X))^{\mathbb{T}})\f} \ar[r] &
{\holim_n((\tilde E\wedge F(E_+,X))^{C_{p^n}})\f} \cr
}$$
where the lower horizontal map is the $p$-completion of the map of the
statement. The top horizontal map is an equivalence by what was shown
above. Moreover, it follows from the proof of the induction step above
that the map of $i$th homotopy groups induced by each of vertical maps
become isomorphisms for $s$ sufficiently large. This uses that
$\pi_tX$ vanishes for $t \ll 0$. This completes the proof.
\end{proof}

Let $G$ be a finite group, and let $T$ be a $G$-spectrum. There is a
natural whole plane spectral sequence called the Tate spectral
sequence
\begin{equation}\label{tatespectralsequence}
\hat{E}_{s,t}^2(G,T) = \hat{H}^{-s}(G,\pi_t(T)) \Rightarrow
\pi_{s+t}(\tate(G,T))
\end{equation}
that converges conditionally in the sense
of~\cite[Def.~5.10]{boardman} to the homotopy groups of the Tate
spectrum and whose $E^2$-term is given by the Tate cohomology of $G$
with coefficients in the homotopy groups of $T$ considered as
$G$-modules. The spectral sequence was first constructed in this
generality in~\cite{greenlees}, but see also~\cite[Sect.~4]{hm2}.

\begin{lemma}\label{colimit}Let $G$ be a finite group and let
$\{T_\alpha\}$ be a filtrered colimit system of $G$-spectra. Suppose
that for all $i\in\Z$, there exists constants $r$ and $c$
{\rm(}independent of $\alpha${\rm)} such that $\hat
E^r_{s,i-s}(G,T_\alpha)=\hat E^{\infty}_{s,i-s}(G,T_\alpha)$, for all
$s\in\Z$, and such that $\hat E^r_{s,i-s}(G,T_\alpha)$ vanishes if
$s<c$. Then the canonical map
$$\hocolim_\alpha\tate(G,T_\alpha) \to
\tate(G,\hocolim_\alpha T_\alpha)$$
is a weak equivalence of spectra.
\end{lemma}

\begin{proof}Let $\Fil_s\tate(G,T_\alpha)$, where $s$ ranges over the
integers, be the filtration of $\tate(G,T_\alpha)$ from
\cite[Sect.~4]{hm2} which gives rise to the Tate spectral
sequence~(\ref{tatespectralsequence}). It induces a filtration of
$\hocolim_\alpha\tate(G,T_\alpha)$, and the canonical map
$$\hocolim_\alpha\tate(G,T_\alpha) \to
\tate(G,\hocolim_\alpha T_\alpha)$$
is filtration preserving. The filtration of the right-hand term gives
rise to the spectral sequence $\hat E^*(G,\hocolim_\alpha T_\alpha)$. 
The filtration of the left-hand term gives rise to the colimit of the
spectral sequences $\hat E^*(G,T_\alpha)$. Since Tate cohomology
preserves filtered colimits, the map of spectral sequences induced by
the map of the statement is an isomorphism. Hence we are done once we
prove that both spectral sequences converge strongly. This is the case
for each of the sequences $\hat E^*(G,T_\alpha)$. Indeed, this follows
from~\cite[Thm.~8.2]{boardman} since the sequences are conditionally
convergent and collapse at the $E^r$-term. Since this is true for all
$\alpha$, the spectral sequence $\hat E^*(G,\hocolim_\alpha T_\alpha)$
also converges strongly. We must show that the sequence
$$\colim_\alpha\hat
E^*(G,T_\alpha)\Rightarrow\pi_*\hocolim_\alpha\tate(G,T_\alpha)$$
converges strongly. It suffices by~\emph{loc.~cit.}~to show that the
spectral sequence converges conditionally which, by definition, means
that
$$\lim_s\colim_\alpha\pi_i\Fil_s\tate(G,T_\alpha)=
R^1\lim_s\colim_\alpha\pi_{i+1}\Fil_s\tate(G,T_\alpha)=0.$$
But the assumptions imply that $\pi_i\Fil_s\tate(G,T_\alpha)$ vanishes
for $s<c$ and for all $\alpha$. The lemma follows.
\end{proof}

We can now prove the following generalization
of~\cite[Addendum~5.4.4]{hm2}.

\begin{prop}\label{trtate}Let $k$ be a perfect field of odd
characteristic $p$, and let $K_0$ be the quotient field of the ring of
Witt vectors $W(k)$. Let $K$ be an algebraic extension of $K_0$ and
let $V$ be the integral closure of $W(k)$ in $K$. Then for all
$n,v \geq 1$ and all $q\geq 0$, the map
$$\hat\Gamma\colon \pi_q(\TR^n(V|K;p),\Z/p^v) \to
\pi_q(\tate(C_{p^n},T(V|K)),\Z/p^v)$$
is an isomorphism.
\end{prop}

\begin{proof}The statement was proved for finite field extensions
in~\emph{loc.~cit.} We write $K=\colim_\alpha K_\alpha$ as the
filtered colimit of all finite subextensions of $K_0$ contained in
$K$.  Let $V_\alpha$ be the integral closure of $W(k)$ in
$K_\alpha$. The colimit of the $V_\alpha$ maps isomorphically onto
$V$, and hence, the induced map of $\mathbb{T}$-spectra
$$\hocolim_\alpha T(V_\alpha|K_\alpha) \to T(V|K)$$
induces an isomorphism on homotopy groups. We claim that also the
induced map of $C_{p^{n-1}}$-fixed point spectra
$$\hocolim_\alpha \TR^n(V_\alpha|K_\alpha;p) \to \TR^n(V|K;p)$$
induces an isomorphism of homotopy groups. To prove the claim, we
recall that by~\cite[Thm.~2.2]{hm}, there is a cofibration sequence of
spectra
$$\borel(C_{p^{n-1}},T(V|K)) \xto{N}
\TR^n(V|K;p) \xto{R}
\TR^{n-1}(V|K;p)$$
that we call the fundamental cofibration sequence. The left-hand term
is the Borel spectrum whose homotopy groups are the abutment of a
natural strongly convergent first-quadrant spectral sequence
$$E_{s,t}^2(C_{p^{n-1}},T(V|K)) = H_s(C_{p^{n-1}},\pi_t(T(V|K)))
\Rightarrow \pi_{s+t}(\borel(C_{p^{n-1}},T(V|K))).$$
It follows that the induced map
$$\hocolim_\alpha \borel(C_{p^n},T(V_\alpha|K_\alpha)) \to
\borel(C_{p^n},\hocolim_\alpha T(V_\alpha|K_\alpha))$$
induces an isomorphism of homotopy groups. The claim then follows by
easy induction over $n \geq 1$.

Finally, we use Lemma~\ref{colimit} to show that the canonical map
$$\hocolim_\alpha\tate(C_{p^n},T(V_\alpha|K_\alpha))\to
\tate(C_{p^n},\hocolim_\alpha T(V_\alpha|K_\alpha))$$
induces an isomorphism of homotopy groups with $\Z/p^v$-coefficients.
It suffices by an easy induction argument to consider the case $v =
1$. Since the Moore spectrum $M_p$ for $\Z/p$ is a finite spectrum, the
canonical map
$$M_p \wedge \tate(C_{p^n},T(V_{\alpha}|K_{\alpha})) \to
\tate(C_{p^n},M_p \wedge T(V_{\alpha}|K_{\alpha}))$$
is a weak equivalence. Hence, it suffices to verify the hypotheses of
Lemma~\ref{colimit} for the colimit system of $C_{p^n}$-spectra
$\{T_{\alpha}\}$, where $T_{\alpha} = M_p \wedge
T(V_{\alpha}|K_{\alpha})$. It suffices as in the proof
of~\cite[Thm.~5.4.3]{hm2} to consider the case where $K_{\alpha}$
contains the $p$th roots of unity. The differential structure of the
spectral sequences
$$\hat{E}_{s,t}^2(C_{p^n},T_{\alpha}) =
\hat{H}^{-s}(C_{p^n},\pi_t(T_{\alpha})) \Rightarrow
\pi_{s+t}(\tate(C_{p^n},T_{\alpha}))$$
was determined completely in~\emph{op.~cit.},~Thm.~5.5.1. In
particular, it was shown there that for all $\alpha$ and for all $r
\geq 2(p^{n+1}-1)/(p-1)$,
$$\hat{E}^r(C_{p^n},T_{\alpha}) = 
\hat{E}^{\infty}(C_{p^n},T_{\alpha}).$$
This verifies the first of the two hypotheses of Lemma~\ref{colimit}.
The structure of the bi-graded abelian group
$\hat{E}^{\infty}(C_{p^n},T_{\alpha})$ was determined
in~\emph{op.~cit.},~Lemma~5.5.3. It follows, in particular, that for
all $\alpha$ and for all $t \geq 2(p^{n+1}-1)/(p-1)-1$,
$$\hat{E}_{s,t}^{\infty}(C_{p^n},T_{\alpha}) = 0.$$
This verfies the second hypothesis of Lemma~\ref{colimit} and hence
the proposition.
\end{proof}

\subsection{}It follows from~\cite[Thm.~3.3.8]{hm2} that $\TR_q^n(\bar
V|\bar{K};p)$ is zero, if $q$ is negative, and that the canonical map
of Witt complexes
$$W_n\Omega_{(\bar{V},\bar{M})}^q \to \TR_q^n(\bar{V}|\bar{K};p)$$
is an isomorphism in degrees $q\leq 2$. Indeed, both sides preserve
filtered colimits. The following result concerns the structure of the
de~Rham-Witt groups on the left-hand side. The definition of the
$W_n(\bar{V})$-module $\smash{{}_hW_n\Omega_{(\bar{V},\bar{M})}^q}$ is
given in~\emph{op.~cit.},~Sect.~3.2.

\begin{prop}\label{drwbar}The group $W_n\Omega_{(\bar{V},\bar{M})}^q$
is divisible, for all $q\geq 1$, and uniquely divisible,
for all $q\geq 2$. Moreover, for all $q\geq 0$, the sequence
is exact:
$$0 \to {}_hW_n\Omega_{(\bar{V},\bar{M})}^q \xto{N}
W_n\Omega_{(\bar{V},\bar{M})}^q \xto{R}
W_{n-1}\Omega_{(\bar{V},\bar{M})}^q \to 0.$$
\end{prop}

\begin{proof}We first show that the sequence of the statement is
exact. It suffices by~\cite[Prop.~3.2.6]{hm2} to show that the
left-hand map is injective. For $q \leq 1$, this follows from the proof
of~\emph{op.~cit.},~Thm.~3.3.8. For $q \geq 2$, we claim that the
canonical map
$$\iota_2 \colon F_*^{n-1}\Omega_{(\bar{V},\bar{M})}^q \to
{}_hW_n\Omega_{(\bar{V},\bar{M})}^q$$
is an isomorphism and that the common group is uniquely
divisible. The injectivity of the map $N$ follows since
the composite
$$F_*^{n-1}\Omega_{(\bar{V},\bar{M})}^q \xto{\iota_2}
{}_hW_n\Omega_{(\bar{V},\bar{M})}^q \xto{N}
W_n\Omega_{(\bar{V},\bar{M})}^q \xto{F^{n-1}}
F_*^{n-1}\Omega_{(\bar{V},\bar{M})}^q$$
is given by multiplication by $p^{n-1}$. To prove the claim, it
suffices to show that the group $\smash{\Omega_{(\bar{V},\bar{M})}^q}$
is divisible, for $q \geq 1$, and uniquely divisible, for $q
\geq 2$. Indeed, this follows immediately
from~\emph{op.~cit.},~Lemma~3.2.5. We write $\bar{K} = 
\colim_{\alpha}K_{\alpha}$ as the filtered colimit of all finite
sub-extension $K_{\alpha}$ of $K$ contained in $\bar{K}$. Then
$\bar{V} = \colim_{\alpha}V_{\alpha}$, where $V_{\alpha}$ is the
integral closure of $V$ in $K_{\alpha}$, and the canonical map
$$\colim_{\alpha} \Omega_{(V_{\alpha},M_{\alpha})}^q \to
\Omega_{(\bar{V},\bar{M})}^q$$
is an isomorphism. It follows from~\emph{op.~cit.},~Lemmas~2.2.3
and~2.2.4 that the common group is divisible, for $q \geq 1$,
and uniquely divisible, for $q \geq 2$, as desired. This completes
the proof that the sequence of the statement is exact.

Since the group $\smash{{}_hW_n\Omega_{(\bar{V},\bar{M})}^q}$ is
uniquely divisible, for $q \geq 2$, an easy induction argument
based on the exact sequence of the statement shows that the group
$\smash{W_n\Omega_{(\bar{V},\bar{M})}^q}$ is uniquely divisible,
for $q \geq 2$, as stated. We will prove in Prop.~\ref{trbar} below
that the group $W_n\Omega_{(\bar{V},\bar{M})}^1$ is divisible.
\end{proof}

We next consider the integral homotopy groups
$$\TH_q(\bar{V}|\bar{K}) = \TR_q^1(\bar{V}|\bar{K};p)$$
of the topological Hochschild spectrum $T(\bar{V}|\bar{K})$.

\begin{lemma}\label{thbar}The group $\TH_q(\bar{V}|\bar{K})$ is
divisible, for all $q>0$, and uniquely divisible, for $q>0$
even. Moreover, for all integers $q$, the canonical map
$$\Omega_{(\bar{V},\bar{M})}^q \to \TH_q(\bar{V}|\bar{K})$$
is a rational isomorophism.
\end{lemma}

\begin{proof}We write $\bar{V} = \colim_{\alpha} V_{\alpha}$ as the
filtered colimit of the integral closure of $V$ in the finite
sub-extensions $K_{\alpha}$ of $K$ contained in $\bar{K}$. Then the
canonical map
$$\colim_{\alpha} \TH_q(V_{\alpha}|K_{\alpha}) \to
\TH_q(\bar{V}|\bar{K})$$
is an isomorphism. We recall from~\cite[Rem.~2.4.2]{hm2} that for $q$
positive and even, the groups $\TH_q(V_{\alpha}|K_{\alpha})$ are
uniquely divisible, and from~\emph{op.~cit.},~Prop.~2.3.4, that
for all integers $q$, the canonical map
$$\Omega_{(V_{\alpha},M_{\alpha})}^q \to
\TH_q(V_{\alpha}|K_{\alpha})$$
is a rational isomorphism. Hence the same holds for the groups
$\TH_q(\bar{V}|\bar{K})$. It remains to show that for $q$ positive and
odd, the group $\TH_q(\bar{V}|\bar{K})$ is divisible. To this end,
we consider the homotopy groups with $\Z/p$-coefficients
$$\TH_q(\bar{V}|\bar{K},\Z/p) = \pi_q(T(\bar{V}|\bar{K}),\Z/p).$$
It follows from~\emph{op.~cit.},~Thm.~B and from the isomorphism
at the beginning of the proof, that there is a canonical isomorphism
of log-differential graded rings
$$\Omega_{(\bar{V},\bar{M})}^*\otimes_{\Z} S_{\Z/p}\{\kappa\}
\xto{\sim}\TH_*(\bar{V}|\bar{K},\Z/p).$$
But, as we saw in the proof of Prop.~\ref{drwbar} above,
$\Omega_{(\bar{V},\bar{M})}^1$ is a divisible group, and hence
this isomorphism becomes
$$S_{\bar{V}/p\bar{V}}\{\kappa\}\xto{\sim}\TH_*(\bar{V}|\bar{K},\Z/p).$$
In particular, the homotopy groups with $\Z/p$-coefficients are
concentrated in even degrees. It follows that for $q$ odd, the middle group in
the short-exact sequence
$$0\to\TH_q(\bar{V}|\bar{K})/p\to\TH_q(\bar{V}|\bar{K},\Z/p)
\to \TH_{q-1}(\bar{V}|\bar{K})[p] \to 0$$
vanishes, and therefore, $\TH_q(\bar{V}|\bar{K})$ is divisible as
stated. Here $A[p]$ denotes the subgroup of elements of order $p$ in
the abelian group $A$.
\end{proof}

\begin{prop}\label{trbar}The group $\TR_q^n(\bar{V}|\bar{K};p)$ is
divisible, if $q\geq 1$, uniquely divisible, if $q\geq 2$
is even, and for all integers $q$, the canonical map
$$W_n\Omega_{(\bar{V},\bar{M})}^q \to
\TR_q^n(\bar{V}|\bar{K};p)$$
is a rational isomorphism. Moreover, the restriction map
$$R\colon \TR_q^n(\bar{V}|\bar{K};p)\to\TR_q^{n-1}(\bar{V}|\bar{K};p)$$
is surjective, for all integers $q$ and all $n \geq 2$.
\end{prop}

\begin{proof}We again use the fundamental cofibration sequence 
$$\borel(C_{p^{n-1}},T(\bar{V}|\bar{K})) \xto{N}
\TR^n(\bar{V}|\bar{K};p) \xto{R}
\TR^{n-1}(\bar{V}|\bar{K};p)$$
which we used above in the proof of Prop.~\ref{trtate}. We shall here
follow the notation from~\cite[Sect.~3.3]{hm2} and write
$${}_h\!\TR_q^n(\bar{V}|\bar{K};p) =
\pi_q(\borel(C_{p^{n-1}},T(\bar{V}|\bar{K}))).$$
The long-exact sequence
$$\cdots\xto{\partial}
{}_h\!\TR_q^n(\bar{V}|\bar{K};p) \xto{N}
\TR_q^n(\bar{V}|\bar{K};p) \xto{R}
\TR_q^{n-1}(\bar{V}|\bar{K};p) \xto{\partial} \cdots$$
and the spectral sequence
$$E^2_{s,t}=H_s(C_{p^{n-1}},F_*^{n-1}\TH_t(\bar{V}|\bar{K}))
\Rightarrow {}_h\!\TR_{s+t}^n(\bar{V}|\bar{K};p)$$
are both sequences of $W_n(\bar{V})$-modules. We recall
from~\emph{op.~cit.},~Lemma~3.3.3 that, in addition, there is a
canonical map of $W_n(\bar{V})$-modules
$${}_hW_n\Omega_{(\bar{V},\bar{M})}^q \to
{}_h\!\TR_q^n(\bar{V}|\bar{K};p)$$
which explains the choice of notation.

We first show that the edge-homomorphism of the spectral sequence
above
$$F_*^{n-1} \TH_q(\bar{V}|\bar{K}) \to
{}_h\!\TR_q^n(\bar{V}|\bar{K};p)$$
is an isomorphism, if $q > 0$ and even, and an injection whose
cokernel is a $p$-torsion group of bounded exponent, if $q > 0$ and
odd. Since the $C_{p^{n-1}}$-action on $T(\bar{V}|\bar{K})$ extends to
an action by the full circle group $\mathbb{T}$, the induced
$C_{p^{n-1}}$-action on homotopy groups is trivial. Hence,
Lemma~\ref{thbar} shows that $E_{s,t}^2$ is zero, if $s > 0$ and $s +
t$ is even, and a $p$-torsion group of bounded exponent, if $s > 0$
and $s + t$ is odd. For degree reasons, the only possible non-zero
differentials are
$$d^r \colon E_{r,q - (r - 1)}^r \to E_{0,q}$$
with $r \geq 2$ and $q$ even. But this is a map from a $p$-torsion
group to a uniquely divisible group and therefore is zero. Hence,
all differentials in the spectral sequence are zero, and the
edge-homomorphism is as stated. It follows that ${}_h\!\TR_q^n(\bar
V|\bar{K};p)$ is a uniquely divisible group, if $q>0$ is even, and
the sum of a uniquely divisible group and a torsion group, if
$q>0$ is odd.

We next show that for all $q$, the canonical map
$${}_hW_n\Omega_{(\bar{V},\bar{M})}^q \to
{}_h\!\TR_q^n(\bar{V}|\bar{K};p)$$
is a rational isomorphism. To this end, we consider the following diagram.
$$\xymatrix{
{F_*^{n-1}\Omega_{(\bar{V},\bar{M})}^q} \ar[r] \ar[d] &
{{}_hW_n\Omega_{(\bar{V},\bar{M})}^q} \ar[d] \cr
{F_*^{n-1}\TH_q(\bar{V}|\bar{K})} \ar[r] &
{{}_h\!\TR_q^n(\bar{V}|\bar{K};p).} \cr
}$$
The lower horizontal map is the edge homomorphism of the
spectral sequence and the remaining maps are the canonical ones. The
top horizontal map is a rational isomorphism by the proof of
Prop.~\ref{drwbar}, and the left-hand vertical map is a rational
isomorphism by Lemma~\ref{thbar}. Finally, we proved above that the
lower horizontal map is a rational isomorphism.

We proceed to show that the restriction map is surjective, or
equivalently, that the long-exact sequence of homotopy groups
associated with the fundamental cofibration sequence breaks into
short-exact sequences
$$0 \to
{}_h\!\TR_q^n(\bar{V}|\bar{K};p) \xto{N} 
\TR_q^n(\bar{V}|\bar{K};p) \xto{R}
\TR_q^{n-1}(\bar{V}|\bar{K};p) \to 0.$$
It suffices to prove that this sequence is exact for $q$ even. Indeed,
this implies that the boundary map in the long-exact sequence is zero,
for all $q$. The sequence above is short-exact for $q = 0$
by~\cite[Thm.~F]{hm}. So we assume that $q \geq 2$ and even and
consider the following diagram.
$$\xymatrix{
{0} \ar[r] &
{{}_hW_n\Omega_{(\bar{V},\bar{M})}^q} \ar[r]^{N} \ar[d] &
{W_n\Omega_{(\bar{V},\bar{M})}^q} \ar[r]^{R} \ar[d] &
{W_{n-1}\Omega_{(\bar{V},\bar{M})}^q} \ar[r] \ar[d] &
{0} \cr
{0} \ar[r] &
{{}_h\!\TR_q^n(\bar{V}|\bar{K};p)} \ar[r]^{N} &
{\TR_q^n(\bar{V}|\bar{K};p)} \ar[r]^{R} &
{\TR_q^{n-1}(\bar{V}|\bar{K};p)} \ar[r] &
{0.}
}$$
The top sequence is exact by Prop.~\ref{drwbar}, and we proved above
that the left-hand vertical map is an isomorphism. Moreover, the lower
left-hand map is injective. Indeed, the composite
$$F_*^{n-1}\TH_q(\bar{V}|\bar{K}) \xto{\sim}
{}_h\!\TR_q^n(\bar{V}|\bar{K};p) \xto{N}
\TR_q^n(\bar{V}|\bar{K};p) \xto{F^{n-1}} 
F_*^{n-1}\TH_q(\bar{V}|\bar{K}),$$
which is given by multiplication by $p^{n-1}$, is an isomorphism, and
we proved above that the left-hand map is an isomorphism. It follows,
by induction on $n \geq 1$, that the lower sequence is exact and that
the middle vertical map is an isomorphism as desired.

It remains to prove that $\TR_q^n(\bar{V}|\bar{K};p)$ is divisible,
if $q>0$ and odd. As in the proof of Lemma~\ref{thbar}, we show that
if $q > 0$ and odd, the middle group in the following short-exact
sequence vanishes.
$$0\to\TR_q^n(\bar{V}|\bar{K};p)/p \to
\TR_q^n(\bar{V}|\bar{K};p,\Z/p) \to
\TR_{q-1}^n(\bar{V}|\bar{K};p)[p] \to 0.$$
To prove this, we recall from Prop.~\ref{trtate} that for $q \geq 0$,
the map
$$\hat{\Gamma} \colon \TR_q^n(\bar{V}|\bar{K};p,\Z/p)
\to \pi_q(\tate(C_{p^n},T(\bar{V}|\bar{K})),\Z/p)$$
is an isomorphism. We use the Tate spectral sequence
$$\hat E^2_{s,t}(C_{p^n}) =
\hat H^{-s}(C_{p^n}, \TH_t(\bar{V}|\bar{K}, \Z/p)) \Rightarrow
\pi_{s+t}(\tate(C_{p^n},\TH(\bar{V}|\bar{K})),\Z/p)$$
to evaluate the right-hand side. As in the proof of Lemma~\ref{thbar},
we find that
$$\hat E^2(C_{p^n}) = \Lambda_{\bar{V}/p\bar{V}}\{u_n\}
\otimes_{\bar{V}/p\bar{V}} S_{\bar{V}/p\bar{V}}\{t^{\pm1},\kappa\}$$
with the generators $u_n$, $t$, and $\kappa$ in bidegrees $(-1,0)$,
$(-2,0)$ and $(0,2)$, respectively. We recall
from~\cite[Sect.~4.4]{hm2} that there is a map
to this spectral sequence from the spectral sequence
$$\hat{E}^2(\mathbb{T}) = S_{\bar{V}/p\bar{V}}\{t^{\pm1},\kappa\}
\Rightarrow \pi_*(\tate(\mathbb{T},T(\bar{V}|\bar{K})),\Z/p)$$
which on the $E^2$-terms is given by the obvious inclusion. For degree
reasons, all differentials in the latter spectral sequence are zero.
It follows that all possible non-zero differentials in the former
spectral sequence are multiplicatively generated from a non-zero
differential on the class $u_n$. Since $T(\bar{V}|\bar{K})$ is a
$T(W(k)|K_0)$-module spectrum,~\cite[Prop.~5.5.4]{hm2} shows that
$$d^{2(\frac{p^{n+1}-1}{p-1})-1}u_n=\mu_n\cdot
(t\kappa)^{\frac{p^{n+1}-1}{p-1}-1}t,$$
where $\mu_n\in\Z/p$ is a unit. We conclude that 
$$\hat E^{\infty}(C_{p^n}) = S_{\bar{V}/p\bar{V}}\{ t^{\pm1}, \kappa\}/
(\kappa^{\frac{p^{n+1}-1}{p-1}-1}),$$
which is concentrated in even degrees as desired.
\end{proof}

\subsection{}We consider the homotopy groups with $\Zp$-coefficients
of $\TR^n(\bar{V}|\bar{K};p)$. We recall that for any spectrum $X$,
the $p$-adic homotopy groups are related to the integral homotopy
groups by a natural short-exact sequence
$$0 \to \Ext(\Qp/\Zp,\pi_q(X)) \to \pi_q(X,\Zp) \to
\Hom(\Qp/\Zp,\pi_{q-1}(X)) \to 0.$$
The right-hand term is also written $T_p(\pi_{q-1}(X))$ and called the
$p$-primary Tate module of the group $\pi_{q-1}(X)$. In the case at
hand, the left-hand term (resp.~the right-hand term) vanishes for $q >
0$ (resp.~for $q > 0$ and odd) since $\TR_q^n(\bar{V}|\bar{K};p)$ is
divisible (resp. torsion-free). In particular, we obtain a canonical
isomorphism
$$\TR_2^n(\bar{V}|\bar{K};p,\Zp) \xto{\sim}
T_p\TR_1^n(\bar{V}|\bar{K};p) \xleftarrow{\sim}
T_pW_n\Omega_{(\bar{V},\bar{M})}^1$$
which, in turn, gives rise to a map of graded
$W_n(\bar{V})\f$-algebras
$$S_{W_n(\bar{V})\f}(T_pW_n\Omega_{(\bar{V},\bar{M})}^1) \to
\TR_*^n(\bar{V}|\bar{K};p,\Zp).$$
We first consider this map for $n = 1$.

\begin{lemma}\label{thcomplete}The $p$-primary Tate module
$T_p\Omega_{(\bar{V},\bar{M})}^1$ is a free $\bar{V}\f\!$-module of
rank one and the canonical map
$$S_{\bar{V}\f}(T_p\Omega_{(\bar{V},\bar{M})}^1) \to
\TH_*(\bar{V}|\bar{K},\Zp)$$
is an isomorphism.
\end{lemma}

\begin{proof}We have a canonical isomorphism
$$T_p\Omega_{(\bar{V},\bar{M})}^1 \xto{\sim} \lim_m
\Omega_{(\bar{V},\bar{M})}^1[p^m]$$
of the $p$-primary Tate module and the limit over $m \geq 1$ of the
$p^m$-torsion subgroups. Let $v^{(m)} \in \bar{V}$, $m \geq 1$, with
$v^{(1)} = -p$ and $(v^{(m+1)})^p = v^{(m)}$. Then
$$\Omega_{(\bar{V},\bar{M})}^1[p^m] = \bar{V}/p^m\bar{V} \cdot
d\log v^{(m)}.$$
Indeed, the case $m = 1$ follows from~\cite[Cor.~2.2.5]{hm2}, and the
general case is proved by induction using the following diagram.
$$\xymatrix{
{ 0 } \ar[r] &
{ \bar{V}/p\bar{V} } \ar[r]^(.46){p^{m-1}} \ar[d]^{d\log v^{(1)}} &
{ \bar{V}/p^m\bar{V} } \ar[r] \ar[d]^{d\log v^{(m)}} &
{ \bar{V}/p^{m-1}\bar{V} } \ar[r] \ar[d]^{d\log v^{(m-1)}} &
{ 0 } \cr
{ 0 } \ar[r] &
{ \Omega_{(\bar{V},\bar{M})}^1[p] } \ar[r] &
{ \Omega_{(\bar{V},\bar{M})}^1[p^m] } \ar[r]^(.46){p} &
{ \Omega_{(\bar{V},\bar{M})}^1[p^{m-1}] } \ar[r] &
{ 0 } \cr
}$$
The lower horizontal sequence is canonically identified with the
coefficients sequence
$$0 \to \TH_2(\bar{V}|\bar{K},\Z/p) \to \TH_2(\bar{V}|\bar{K},\Z/p^m)
\to \TH_2(\bar{V}|\bar{K},\Z/p^{m-1}) \to 0$$
which is exact since $\TH_q(\bar{V}|\bar{K},\Z/p^v)$ is zero for $q$
odd. It follows that the $p$-primary Tate module
$T_p\Omega_{(\bar{V},\bar{M})}^1$ is a free $\bar{V}\f$-module of rank
one and that $\kappa=\{d\log v^{(m)}\}_{m\geq 1}$ is a
generator. It remains to show that the composite map
$$S^2\wedge\TH(\bar{V}|\bar{K})\f\xto{\kappa\wedge\id}
\TH(\bar{V}|\bar{K})\f\wedge\TH(\bar{V}|\bar{K})\f\xto{\mu}
\TH(\bar{V}|\bar{K})\f$$
induces an isomorphism of homotopy groups in degrees $q \geq 2$. Since
this is a map between $p$-complete spectra, it suffices to show that
the induced map of homotopy groups with $\Z/p$-coefficients is an
isomorphism in degrees $q \geq 2$. But this follows
from~\cite[Thm.~B]{hm2}. 
\end{proof}

The following result together with Prop.~\ref{trbar} above completes
the proof of Thm.~\ref{gradedringstructure} of the induction.

\begin{prop}\label{trcomplete}The $p$-primary Tate module
$\smash{T_pW_n\Omega_{(\bar{V},\bar{M})}^1}$ is a free module of rank
one over $\smash{W_n(\bar{V})\f}$ and the canonical map
$$S_{W_n(\bar{V})\f}(T_pW_n\Omega_{(\bar{V},\bar{M})}^1) \to
\TR_*^n(\bar{V}|\bar{K};p,\Zp)$$
is an isomorphism.
\end{prop}

\begin{proof}We know from Prop.~\ref{trtate} that the map
$$\hat{\Gamma} \colon \TR^n(\bar{V}|\bar{K};p)
\to \tate(C_{p^n},\TH(\bar{V}|\bar{K}))$$
induces an isomorphism of homotopy groups with $\Zp$-coefficients in
non-negative degrees. It follows from Lemma~\ref{thcomplete} that
the Tate spectral sequence
$$\hat E^2_{s,t}=\hat H^{-s}(C_{p^n},\TH_t(\bar{V}|\bar{K},\Zp))
\Rightarrow\pi_{s+t}(\tate(C_{p^n},\TH(\bar{V}|\bar{K})),\Zp)$$
takes the form $\hat E^2 = S_{\bar{V}/p^n\bar{V}}\{t^{\pm 1},\kappa\}$
with $t$ and $\kappa$ in bi-degrees $(-2,0)$ and $(0,2)$,
respectively. Since all non-zero elements are located in even total
degree, all differentials are zero and the groups
$\TR_q^n(\bar{V}|\bar{K};p,\Zp)$ are concentrated in even degrees.
Moreover, the spectral sequence is multiplicative, and
multiplication by $t^{-1}$ induces an isomorphism of $\hat
E_{s,t}^{\infty}$ onto $\hat E_{s+2,t}^{\infty}$. Let
$\alpha_n\in\TR_2^n(\bar{V}|\bar{K};p,\Zp)$ be a homotopy 
class such that $\smash{\hat\Gamma(\alpha_n)}$ is represented in the
spectral sequence by the element $t^{-1}$. Then multiplication by
$\alpha_n$ induces an isomorphism
$$\TR_q^n(\bar{V}|\bar{K};p,\Zp) \xto{\sim}
\TR_{q+2}^n(\bar{V}|\bar{K};p,\Zp)$$
for all positive integers $q$. This proves the proposition.
\end{proof}

\begin{addendum}\label{frobeniussurjective}The Frobenius
$$F\colon \TR_q^n(\bar{V}|\bar{K};p,\Zp)\to\TR_q^{n-1}(\bar{V}|\bar{K};p,\Zp)$$
is surjective.
\end{addendum}

\begin{proof}For $q=0$, the statement is that $F\colon W_n(\bar{V})\f\to
W_{n-1}(\bar{V})\f$ is surjective which we proved in
Cor.~\ref{frobsurj}. Hence it suffices by Prop.~\ref{trcomplete} to
show that
$$F \colon \TR_2^n(\bar{V}|\bar{K};p,\Zp) \to
\TR_2^{n-1}(\bar{V}|\bar{K};p,\Zp)$$
takes a $W_n(\bar{V})\f$-module generator to a
$W_{n-1}(\bar{V})\f$-module generator. By Prop.~\ref{trtate} we may
instead show that the map
$$F \colon \pi_2(\tate(C_{p^n},T(\bar{V}|\bar{K})),\Zp) \to
\pi_2(\tate(C_{p^{n-1}},T(\bar{V}|\bar{K})),\Zp)$$
takes a $W_n(\bar{V})\f$-module generator to a
$W_{n-1}(\bar{V})\f$-module generator. It follows from the
construction of the Tate spectral sequence that the latter map induces
a map of Tate spectral sequences that on $E^2$-terms is given by the map
$$S_{\bar{V}/p^n\bar{V}}\{t^{\pm 1},\kappa\} \to
S_{\bar{V}/p^{n-1}\bar{V}}\{t^{\pm 1},\kappa\}$$
that takes generators $t$ and $\kappa$ on the left-hand side to the
generators $t$ and $\kappa$ on the right-hand side and that maps
$\bar{V}/p^n\bar{V}$ to $\bar{V}/p^{n-1}\bar{V}$ by the canonical
projection; see~\cite[Sect.~4]{hm2}. The element $t^{-1}$ on the
left-hand side represents a $W_n(\bar{V})\f$-module generator of 
$\pi_2(\tate(C_{p^n},T(\bar{V}|\bar{K})),\Zp)$ and the element
$t^{-1}$ on the right-hand side represents a
$W_{n-1}(\bar{V})\f$-module generator of
$\smash{\pi_2(\tate(C_{p^{n-1}},T(\bar{V}|\bar{K})),\Zp)}$.
This completes the proof.
\end{proof}

%
%
%
%

We conclude this paragraph with a brief discussion of the
homotopy groups with $\Qp/\Zp$-coefficients of the spectrum
$\TR^n(\bar{V}|\bar{K};p)$. 

\begin{lemma}\label{tatemodule}Let $A$ be a divisible abelian
group and let $T_p(A)$ be the $p$-primary Tate module of $A$. Then
there is a canonical isomorphism
$$\Tor(A,\Qp/\Zp) \xto{\sim} T_p(A) \otimes \Qp/\Zp.$$
\end{lemma}

\begin{proof}Since $A$ is a divisible abelian group, the short-exact
sequence
$$0 \to \Z \to \Z[\textstyle{\frac{1}{p}}] \to \Qp/\Zp \to 0$$
induces a short-exact
$$0 \to \Hom(\Qp/\Zp,A) \to \Hom(\Z[\textstyle{\frac{1}{p}}],A) \to
\Hom(\Z,A) \to 0$$
that is usually written
$$0 \to T_p(A) \to V_p(A) \to A \to 0.$$
This sequence, in turn, gives rise to a six-term exact sequence 
$$\begin{aligned}
0 & \to \Tor(T_p(A),\Qp/\Zp) \to
\Tor(V_p(A),\Qp/\Zp) \to
\Tor(A,\Qp/\Zp) \cr
{} & \to T_p(A) \otimes \Qp/\Zp \to
V_p(A) \otimes \Qp/\Zp \to
A \otimes \Qp/\Zp \to 0. \cr
\end{aligned}$$
To prove the statement of the lemma, it suffices to prove that the
second and fifth terms of this sequence vanish. But this follows
immediately from the exact sequence at the beginning of the proof
considered as a flat resolution of $\Qp/\Zp$ and from the fact that
$V_p(A)$ is a $\Z[\frac{1}{p}]$-module.
\end{proof}

\begin{cor}For all integers $q$ and all positive integers $n$, there
is a canonical isomorphism of abelian groups
$$\TR_q^n(\bar{V}|\bar{K};p,\Qp/\Zp) \xto{\sim}
\TR_q^n(\bar{V}|\bar{K};p,\Zp) \otimes \Qp/\Zp.$$
\end{cor}

\begin{proof}This follows immediately from
Thm.~\ref{gradedringstructure}, Lemma~\ref{tatemodule}, and from the
definition of homotopy groups with $\Qp/\Zp$-coefficients which we
recalled at the beginning of Sect.~\ref{algclosed} above.
\end{proof}

\subsection{}Let $\e=\{\e^{(v)}\}_{v\geq 1}$ be a sequence of
primitive $p^{v-1}$th roots of unity in $\bar{K}$ that are compatible
in the sense that $(\e^{(v+1)})^p=\e^{(v)}$. The sequence $\e$ is a
generator of the $p$-primary Tate module $T_p(\bar{K}^*)$, and the
image by the canonical isomorphism
$$T_p(\bar{K}^*) \xto{\sim} T_pK_1(\bar{K})$$
is the associated $K$-theory Bott element $\beta_{\e} =
\beta_{\e}^K$.

\begin{lemma}The image of the $K$-theory Bott element $\beta_\e^K$ by
the map
$$T_pK_1(\bar{K}) \to T_p\TR_1^n(\bar{V}|\bar{K};p)$$
induced by the cyclotomic trace is the element
$\beta_{\e,n} = \beta_{\e,n}^{\TR}
= \{ d\log_n\e^{(v)} \}_{v \geq 1}$.
\end{lemma}

\begin{proof}We recall from~\cite[Sect.~2.3]{hm2} that the map
$$d\log_n \colon \bar{M} = \bar{V} \cap \bar{K}^* \to
\TR_1^n(\bar{V}|\bar{K};p)$$
is defined to be the composite of the canonical inclusion
$\alpha \colon \bar{M} \hookrightarrow \bar{K}^*$, the canonical
isomorphism $\bar{K}^* \to K_1(\bar{K})$, and the cyclotomic
trace $K_1(\bar{K}) \to \TR_1^n(\bar{V}|\bar{K};p)$. The stated
formula is now clear.
\end{proof}

The $\TR$-theory Bott element
$$\beta_{\e,n} = \beta_{\e,n}^{\TR} \in T_p\TR_1^n(\bar{V}|\bar{K};p)
\xleftarrow{\sim} T_pW_n\Omega_{(\bar{V},\bar{M})}^1$$
satisfies that $R(\beta_{n,\e})=F(\beta_{n,\e})=\beta_{n-1,\e}$ and
hence defines a Bott element
$$\beta_{\e} = \beta_{\e}^{\TF} = \{ \beta_{\e,n} \}_{n\geq 1}
\in\lim_FT_pW_n\Omega_{(\bar{V},\bar{M})}^1.$$
The right-hand group is a free module of rank one over the ring
$\lim_F W_n(\bar{V})\f$ which we identify with the ring
$W(R_{\bar{V}})$ via ring isomorphism
$$\psi \colon W(R_{\bar{V}}) \xto{\sim}
\lim_F W_n(\bar{V})\f$$
of Addendum~\ref{thelimitwittring}. The Bott element $\beta_{\e}$,
however, is not a generator. Instead we shall prove the following result.
Let $\e_n = \{\e^{(v+n)}\}_{v\geq 1}$. We denote by $\e_n$ the
corresponding element of the ring $R_{\bar{V}}$ and by $[\e_n]$ the
associated Teichm\"uller representative of the Witt ring
$W(R_{\bar{V}})$. 

\begin{prop}\label{generator}There exists a unique
$W(R_{\bar{V}})$-module generator
$$\alpha_{\e} \in \lim_F T_p W_n\Omega_{(\bar{V},\bar{M})}^1$$
with the property that $\beta_{\e}=\psi([\e_1]-1) \cdot \alpha_{\e}$.
\end{prop}

The uniqueness part of Prop.~\ref{generator} follows immediately from
Addendum~\ref{domain}; the existence part is proved in
Sect.~\ref{generatorsection} below. We conclude this section with the
proofs of Thm.~\ref{wittcomplexstructure} and~\ref{tcthm} of the
introduction.

\begin{proof}[Proof of Thm.~\ref{wittcomplexstructure}]We first prove
the statement (i) of Thm.~\ref{wittcomplexstructure}. It follows from
Addendum~\ref{frobeniussurjective} that the canonical projection 
$$\pr_{n,1} \colon \lim_F T_pW_n\Omega_{(\bar{V},\bar{M})}^1 \to
T_pW_n\Omega_{(\bar{V},\bar{M})}^1$$
is surjective. Similarly, Cor.~\ref{frobsurj} shows that the canonical
projection
$$\pr_{n,0} \colon \lim_F W_n(\bar{V})\f \to W_n(\bar{V})\f$$
is surjective. Moreover, the domain (resp.~the range) of the map
$\pr_{n,1}$ is a free module of rank one over the domain (resp.~range)
of the map $\pr_{n,0}$. It follows that the map $\pr_{n,1}$ takes a
generator to a generator. Hence, the class defined by
$$\alpha_{\e,n} = \pr_{n,1}(\alpha_{\e}) \in
T_p W_n\Omega_{(\bar{V},\bar{M})}^1$$
is a $W_n(\bar{V})\f$-module generator. Finally,
$$\begin{aligned}
\beta_{\e,n} & = \pr_{n,1}(\beta_{\e})
= \pr_{n,1}(\psi([\e_1]-1) \cdot \alpha_{\e}) \cr
{} & = \pr_{n,0}(\psi([\e_1] - 1)) \cdot \pr_{n,1}(\alpha_{\e})
= \theta_n([\e_n] - 1) \cdot \alpha_{\e,n} \cr
\end{aligned}$$
as stated.

We next consider the statement~(ii) of
Thm.~\ref{wittcomplexstructure}. It is immediate from the definition
of $\alpha_{\e,n}$ that $F(\alpha_{\e,n}) = \alpha_{\e,n-1}$. The
restriction map induces a self-map
$$R \colon \lim_F T_pW_n\Omega_{(\bar{V},\bar{M})}^1
\to \lim_F T_pW_n\Omega_{(\bar{V},\bar{M})}^1.$$
of a free $W(R_{\bar{V}})$-module of rank one. Hence there exists an
element $\xi \in W(R_{\bar{V}})$ such that
$R(\alpha_{\e}) = \psi(\xi) \cdot \alpha_{\e}$, and since $W(R_{\bar
V})$ is an integral domain, this element $\xi$ is unique. Since
$R(\beta_{\e}) = \beta_{\e}$ and $\beta_{\e} = \psi([\e_1] - 1) \cdot
\alpha_{\e}$, we find that
$$\psi([\e_1] - 1) \cdot \alpha_{\e} = 
R(\psi([\e_1] - 1) \cdot \alpha_{\e}) = 
\psi([\e_2] - 1) \cdot R(\alpha_{\e})$$
which shows that $\xi=([\e_1]-1)/([\e_2]-1)$. The stated formula now
follows from Addendum~\ref{thelimitwittring}.

The proof of statement~(iii) of Thm.~\ref{wittcomplexstructure} is
similar to the proof of statement~(ii) above. We have
$\beta_{\e}^{\sigma} = \chi(\sigma) \cdot \beta_{\e}$ and $\beta_{\e}
= \psi([\e_1]-1) \cdot \alpha_{\e}$ which shows that
$$\psi([\e_1^{\sigma}] - 1) \cdot \alpha_{\e}^{\sigma}
= \chi(\sigma) \cdot \psi([\e_1] - 1) \cdot \alpha_{\e}.$$
This is an equation among elements of a free module of rank one over
an integral domain. Hence, we obtain
$$\alpha_{\e}^{\sigma} = \chi(\sigma) \cdot \psi(([\e_1] -
1)/([\e_1^{\sigma}] - 1)) \cdot \alpha_{\e}$$
and apply Addendum~\ref{thelimitwittring} to obtain the stated formula.
\end{proof}

\begin{proof}[Proof of Thm.~\ref{tcthm}]Let $\TF(\bar{V}|\bar{K};p)$
be the homotopy limit spectrum
$$\TF(\bar{V}|\bar{K};p) = \holim_{F} \TR^n(\bar{V}|\bar{K};p)$$
where the structure map in the limit system is the Frobenius map. The
restriction and Frobenius maps give rise to self-maps of the spectrum
$\TF(\bar{V}|\bar{K};p)$ and the topological cyclic homology spectrum
$\TC(\bar{V}|\bar{K};p)$ is defined to be the homotopy
equalizer. Hence, the homotopy groups form a long-exact sequence
$$\ldots \to \TC_q(\bar{V}|\bar{K};p) \to
\TF_q(\bar{V}|\bar{K};p) \xto{1 - R}
\TF_q(\bar{V}|\bar{K};p) \xto{\partial}
\TC_{q-1}(\bar{V}|\bar{K};p) \to \ldots$$
and we have the Milnor short-exact sequence
$$0 \to R^1\lim_F \TR_{q+1}^n(\bar{V}|\bar{K};p) \to
\TF_q(\bar{V}|\bar{K};p) \to
\lim_F \TR_q^n(\bar{V}|\bar{K};p) \to 0.$$
The same is true for the homotopy groups with $\Zp$-coefficients. It
follows from Thms.~\ref{gradedringstructure}
and~\ref{wittcomplexstructure} that the map
$$g_{\e,m} \colon W(R_{\bar{V}}) \to \TF_{2m}(\bar{V}|\bar{K};p,\Zp)$$
that takes $a$ to $\psi(a) \cdot \alpha_{\e}^m$ is an isomorphism,
for all non-negative integers $m$, and that the groups in odd degrees
are zero. Moreover, the following diagram commutes.
$$\xymatrix{
{ W(R_{\bar{V}}) } \ar[rr]^(.37){g_{\e,m}}
\ar[d]^{(\frac{[\e_1]-1}{[\e_2]-1})^m W(\varphi^{-1})} &
{} &
{ \TF_{2m}(\bar{V}|\bar{K};p,\Zp) } \ar[d]^{R} \cr
{ W(R_{\bar{V}}) } \ar[rr]^(.37){g_{\e,m}} &
{} &
{ \TF_{2m}(\bar{V}|\bar{K};p,\Zp). } \cr
}$$
The statement of the theorem now follows from Cor.~\ref{calctc}.
\end{proof}

\section{The generator $\alpha_{\e}$}\label{generatorsection}

\subsection{}We obtain the generator $\alpha_{\e}$ from the Thom class
$\lambda_{-L}$ in topological $K$-theory associated with negative
the canonical line bundle over $\mathbb{P}^{\infty}(\C)$. We first
recall some standard facts about topological $K$-theory.

Let $ku$ be the connective cover of the $p$-completion of the spectrum
representing periodic complex $K$-theory. This is a ring spectrum
whose homotopy groups
$$\pi_*(ku) = S_{\Zp}\{ \beta \}$$
is the polynomial algebra over $\Zp$ on the Bott element
$\beta = \beta^{ku}$ of degree $2$.

The canonical projection $f \colon \mathbb{T} \to \{1\}$ from the
circle group to the trivial group gives rise to an adjoint pair of
functors $(f^*,f_*)$ between the stable homotopy category and the
$\mathbb{T}$-stable homotopy category. (The functor $f_*$ was left out
of the notation in the bottom row of the diagram at the beginning of
the introduction.) We consider the Tate spectrum of the
$\mathbb{T}$-spectrum $f^*ku$ obtained from the spectrum $ku$ in this
way. The Tate spectral sequence takes the form
$$\hat{E}^2 = S_{\Zp}\{ t^{\pm1}, \beta \} \Rightarrow
\pi_*(\tate(\mathbb{T},f^*ku))$$
with the generators $t$ and $\beta$ located in bi-degrees $(-2,0)$ and
$(0,2)$, respectively. Since the non-zero groups are all located in
even total degrees, it follows that all differentials are zero. Hence,
the graded ring $\pi_*(\tate(\mathbb{T},f^*ku))$ is the tensor product
of a power series algebra and a Laurent polynomial algebra over $\Zp$
on generators $h$ and $\lambda_{-L}$ that represent $t\beta$ and
$t^{-1}$, respectively. We now explicitly choose a pair of generators
$h$ and $\lambda_{-L}$.

Let $L = \mathcal{O}(1)$ be the canonical line bundle over
$\mathbb{P}^{\infty}(\C)$ and define $\mathbb{P}_{-s}^{\infty}(\C)$ to
be the Thom spectrum of the virtual bundle $-sL$. We recall that
$\mathbb{P}_{-s}^{\infty}(\C)$ may be given the structure of a
CW-spectrum with one cell in every even dimension greater than or
equal to $-2s$. Moreover, it follows
from~\cite[Prop.~II.4.4]{lewismaysteinberger} that for all integers
$s$, there is a canonical isomorphism
$$\pi_{2s}(\tate(\mathbb{T},f^*ku)) \xto{\sim}
K^{-2s}(\mathbb{P}_{-s}^{\infty}(\C),\Zp).$$
The right-hand side is a free
$K^0(\mathbb{P}_0^{\infty}(\C),\Zp)$-module of rank one generated by
the standard Thom class $\lambda_{-sL}$ of the virtual bundle
$-sL$. Moreover, the product
$$K^{-2s}(\mathbb{P}_{-s}^{\infty}(\C),\Zp) \otimes
K^{-2s'}(\mathbb{P}_{-s'}^{\infty}(\C),\Zp) \to
K^{-2(s+s')}(\mathbb{P}_{-(s+s')}^{\infty}(\C),\Zp)$$
takes $\lambda_{-sL} \otimes \lambda_{-s'L}$ to
$\lambda_{-(s+s')L}$. This defines the generator $\lambda_{-L}$. We
define the generator $h = [L] - 1 \in
K^0(\mathbb{P}_0^{\infty}(\C),\Zp)$ to be the class of the reduced 
canonical line bundle. Let $\hat{\beta} \in
K^{-2}(\mathbb{P}_{-1}^{\infty}(\C),\Zp)$ be the image of the 
Bott element $\beta$ by the composite of the canonical maps
$$K^{-2}(S^0,\Zp) \xto{\gamma}
K^{-2}(\mathbb{P}_0^{\infty}(\C),\Zp) \to 
K^{-2}(\mathbb{P}_{-1}^{\infty}(\C),\Zp).$$
The class $\hat{\beta}$ is represented in the spectral sequence by the
element $\beta$ and the product
$$K^0(\mathbb{P}_0^{\infty}(\C),\Zp) \otimes
K^{-2}(\mathbb{P}_{-1}^{\infty}(\C),\Zp) \to
K^{-2}(\mathbb{P}_{-1}^{\infty}(\C),\Zp)$$
takes $([L] - 1) \otimes \lambda_{-L}$ to $\hat{\beta}$.

\begin{prop}\label{topologicalktheory}The graded ring
$$\pi_*(\tate(\mathbb{T},f^*ku)) = \bigoplus_{s \in \Z}
K^{-2s}(\mathbb{P}_{-s}^{\infty}(\C),\Zp)$$
is the tensor product of the power series algebra and the
Laurent polynomial algebra over $\Zp$ generated by the class
$[L] - 1 \in K^0(\mathbb{P}_0^{\infty}(\C),\Zp)$ of the reduced
canonical line bundle and by the standard Thom class $\lambda_{-L} \in
K^{-2}(\mathbb{P}_{-1}^{\infty}(\C),\Zp)$ of negative the canonical
line bundle, respectively. Moreover
$$\hat{\beta} = ([L] - 1) \cdot \lambda_{-L}.$$
\end{prop}

\subsection{}It follows from theorems of Suslin~\cite{suslin,suslin1}
that the spectrum $ku$ and the $p$-completion of the algebraic
$K$-theory spectrum $K(\bar{K})$ are weakly equivalent ring spectra
and that we can choose a weak equivalence of ring spectra
$$\iota_{\e} \colon ku \xto{\sim} K(\bar{K})\f$$
such that the induced map of homotopy groups takes the Bott element
$\beta$ to the Bott element $\beta_{\e}$. We obtain a map of
$\mathbb{T}$-ring spectra
$$\tau_{\e} \colon f^*ku \to T(\bar{V}|\bar{K})\f$$
as the adjoint of the composite of the map $\iota_{\e}$ and the trace map
$$\tr \colon K(\bar{K})\f \to f_*T(\bar{V}|\bar{K})\f.$$
We also recall from Lemma~\ref{Ttate} that the canonical map
$$\tate(\mathbb{T},T(\bar{V}|\bar{K}))\f \to
\holim_F \tate(C_{p^n},T(\bar{V}|\bar{K}))\f$$
is a weak equivalence and from Prop.~\ref{trtate} that the map
$$\hat{\Gamma}_n \colon \TR^n(\bar{V}|\bar{K};p)\f \to
\tate(C_{p^n},T(\bar{V}|\bar{K}))\f$$
induces an isomorphism of homotopy groups in non-negative
degrees. Moreover, the definition of the map
$\hat{\Gamma}_n$~\cite[Sect.~1.1]{hm2} shows that $F \circ
\hat{\Gamma}_n = \hat{\Gamma}_{n-1} \circ F$.  Hence there is a
well-defined multiplicative map in the stable homotopy category
$$\hat{\Gamma} \colon \TF(\bar{V}|\bar{K};p)\f = \holim_F
\TR^n(\bar{V}|\bar{K};p)\f \to 
\tate(\mathbb{T},T(\bar{V}|\bar{K}))\f$$
which induces an isomorphism of homotopy groups in non-negative
degrees. Let
$$\hat{\tau}_{\e} \colon \tate(\mathbb{T},f^*ku) \to
\tate(\mathbb{T},T(\bar{V},\bar{K})\f)$$
be the map of $\mathbb{T}$-Tate spectra induced from the map
$\tau_{\e}$.

\begin{lemma}\label{tauhatbeta}The map $\hat{\tau}_{\e}$ takes the
class $\hat{\beta}$ to the class $\hat{\Gamma}(\beta_{\e})$.
\end{lemma}

\begin{proof}We recall~\cite[Sect.~1.1]{hm2} the following commutative
diagram of spectra and multiplicative maps.
$$\xymatrix{
{ \TR^n(\bar{V}|\bar{K};p) } \ar[r]^{R} \ar[d]^{\Gamma_n} &
{ \TR^{n-1}(\bar{V}|\bar{K};p) } \ar[d]^{\hat{\Gamma}_{n-1}} \cr
{ \coborel(C_{p^{n-1}}, T(\bar{V}|\bar{K})) } \ar[r]^{R^h} &
{ \tate(C_{p^{n-1}}, T(\bar{V}|\bar{K})). } \cr
}$$
All maps in this diagram are compatible with the respective Frobenius
maps. Hence, we obtain a well-defined commutative diagram of
spectra and multiplicative maps in the stable homotopy category
\begin{equation}\label{RRh}
\xymatrix{
{ \TF(\bar{V}|\bar{K};p)\f } \ar[r]^{R} \ar[d]^{\Gamma} &
{ \TF(\bar{V}|\bar{K};p)\f } \ar[d]^{\hat{\Gamma}} \cr
{ \coborel(\mathbb{T},T(\bar{V}|\bar{K}))\f } \ar[r]^{R^h} &
{ \tate(\mathbb{T},T(\bar{V}|\bar{K}))\f } \cr
}
\end{equation}
and Prop.~\ref{trtate} shows that the vertical maps induce
isomorphisms of homotopy groups in non-negative degrees. The Bott
element $\hat{\beta}$ is defined to be the image of the Bott element
$\beta = \beta^{ku}$ by the map of homotopy groups induced by the
following composite map
$$ku \xto{\gamma} \coborel(\mathbb{T},f^*ku) \xto{R^h}
\tate(\mathbb{T},f^*ku).$$
Moreover, the following diagram commutes
$$\xymatrix{
{ ku } \ar[r]^(.4){\gamma} \ar[d]^{\tr \circ \iota_{\e}} &
{ \coborel(\mathbb{T},f^*ku) } \ar[r]^{R^h} \ar[d]^{\tau_{\e}'} &
{ \tate(\mathbb{T},f^*ku) } \ar[d]^{\hat{\tau}_{\e}} \cr
{ \TF(\bar{V}|\bar{K};p)\f } \ar[r]^(.45){\Gamma} &
{ \coborel(\mathbb{T},T(\bar{V}|\bar{K})\f) } \ar[r]^{R^h} &
{ \tate(\mathbb{T},T(\bar{V}|\bar{K})\f) } \cr
}$$
The map of homotopy groups induced by the left-hand vertical map takes
the Bott element $\beta = \beta^{ku}$ to the Bott element $\beta_{\e}
= \beta_{\e}^{\TF}$. Since the restriction map
$$R \colon \TF_2(\bar{V}|\bar{K};p,\Zp) \to
\TF_2(\bar{V}|\bar{K};p,\Zp)$$
fixes $\beta_{\e}$, and since the
diagram~(\ref{RRh}) commutes, the composition of the lower horizontal
maps in the diagram above takes the Bott element $\beta_{\e}$ to the
class $\hat{\Gamma}(\beta_{\e})$. The lemma follows.
\end{proof}

\begin{lemma}\label{tauhatL}The map $\hat{\tau}_{\e}$ takes the class
$[L]$ to the class $\hat{\Gamma}(\psi([\e_1']))$ for a possibly
different sequence $\e' = \{ \e'{}^{(n)} \}_{n \geq 1}$ of primitive
$p^{n-1}$th roots of unity in $\bar{V}\f$.
\end{lemma}

\begin{proof}We recall from Addendum~\ref{thelimitwittring} that the
restriction map
$$R \colon \TF_0(\bar{V}|\bar{K};p,\Zp) \to
\TF_0(\bar{V}|\bar{K};p,\Zp)$$
takes $\psi([\e'])$ to $\psi([\e_1'])$. Hence, in view of the
diagram~(\ref{RRh}), it suffices to prove that the map of homotopy
groups induced by the map
$$\tau_{\e}' \colon \coborel(\mathbb{T},f^*ku) \to
\coborel(\mathbb{T},T(\bar{V}|\bar{K})\f)$$
takes the class $[L]$ to $\Gamma(\psi([\e']))$. Let
$$\tau_{\e,n}' \colon \coborel(C_{p^{n-1}},f^*ku) \to
\coborel(C_{p^{n-1}},T(\bar{V}|\bar{K})\f)$$
be the map induced by the map $\tau_{\e}$. The maps $\tau_{\e,n}'$ are
compatible with the Frobenius maps such that we have a map of
pro-spectra
$$\{ \tau_{\e,n}' \}_{n \geq 1} \colon
\{ \coborel(C_{p^{n-1}},f^*ku) \}_{n \geq 1} \to
\{ \coborel(C_{p^{n-1}},T(\bar{V}|\bar{K})\f) \}_{n \geq 1}.$$
Moreover, there is a commutative diagram of spectra
$$\xymatrix{
{ \coborel(\mathbb{T},f^*ku) } \ar[r]^(.45){\tau_{\e}'} \ar[d] &
{ \coborel(\mathbb{T},T(\bar{V}|\bar{K})\f) } \ar[d] \cr
{ \holim_F \coborel(C_{p^{n-1}},f^*ku) } \ar[r] &
{ \holim_F \coborel(C_{p^{n-1}},T(\bar{V}|\bar{K})\f) } \cr
}$$
where the lower horizontal map is induced by the map of pro-spectra
$\{\tau_{\e,n}'\}$. The vertical maps are weak equivalences by the
proof of Lemma~\ref{Ttate}. The Atiyah-Segal completion theorem shows
that as a ring
$$\pi_0(\coborel(C_{p^{n-1}},f^*ku)) = \Zp[t_n]/(t_n^{p^{n-1}} - 1)$$
where $t_n$ is the image of $[L]$, and Prop.~\ref{trtate}
shows that the ring homomorphism
$$\Gamma_n \colon W_n(\bar{V})\f \to
\pi_0(\coborel(C_{p^{n-1}},T(\bar{V}|\bar{K})\f))$$
is an isomorphism. Hence the map of pro-spectra $\{\tau_{\e,n}'\}$
induces a map of pro-rings
$$\tilde{a} = \{ \tilde{a}^{(n)} \}_{n \geq 1} \colon
\{ \Zp[t_n]/(t_n^{p^{n-1}}) \}_{n \geq  1} \to 
\{ W_n(\bar{V})\f \}_{n \geq 1}.$$
The map $\tilde{a}$ determines and is determined by the sequence $a
= \{ a^{(n)} \}_{n \geq 1}$, where $a^{(n)} = \tilde{a}^{(n)}(t_n)$ is
a $p^{n-1}$th root of unity in $W_n(\bar{V})\f$ and where $F(a^{(n)})
= a^{(n-1)}$. We  claim that for all $n \geq 1$, the map
$\tilde{a}^{(n)}$ is injective or equivalently that $a^{(n)}$ is a
primitive $p^{n-1}$th root of unity. Indeed, the map $\tau_{\e,n}'$
induces a map from the spectral sequence
$$E_{s,t}^2 = H^{-s}(C_{p^{n-1}},\pi_t(ku)) \Rightarrow
\pi_{s+t}(\coborel(C_{p^{n-1}},f^*ku))$$
to the spectral sequence
$$E_{s,t}^2 = H^{-s}(C_{p^{n-1}},T(\bar{V}|\bar{K})\f) \Rightarrow
\pi_{s+t}(\coborel(C_{p^{n-1}},T(\bar{V}|\bar{K})\f)).$$
The map of $E^2$-terms is injective and, for degree reasons, all
differentials in the two spectral sequences are zero, so the induced
map of $E^{\infty}$-terms is also injective. Hence the induced map of
homotopy groups is injective as claimed. Finally,
Lemma~\ref{rootsofunity} shows that $a = \{ [\e'{}^{(n)}]_n\}_{n \geq 1}$,
where $\e' = \{ \e'{}^{(n)} \}_{n   \geq 1}$ is a sequence of
primitive $p^{n-1}$th roots of unity in $\bar{V}\f$ such that
$(\e'{}^{(n)})^p = \e'{}^{(n-1)}$. The statement follows since the
isomorphism $\psi$ of Addendum~\ref{thelimitwittring} takes $[\e']$ to
$\{ [\e'{}^{(n)}]_n \}_{n \geq 1}$.
\end{proof}

\begin{remark}We expect that the two sequences $\e'$ and $\e$ in the
statement of Lemma~\ref{tauhatL} are equal.
\end{remark}

\begin{proof}[Proof of Prop.~\ref{generator}]We recall from the proof
  of Prop.~\ref{trcomplete} that the Tate spectral sequence associated
  with the $\mathbb{T}$-spectrum $T(\bar{V}|\bar{K})$ takes the form
$$\hat{E}^2 = S_{\bar{V}\f}\{t^{\pm 1},\kappa\}
\Rightarrow \pi_*(\tate(\mathbb{T},T(\bar{V}|\bar{K})),\Zp)$$
where the generators $t$ and $\kappa$ are located in bi-degrees
$(-2,0)$ and $(2,0)$. Since the non-zero groups are concentrated in
even total degree, all differentials are zero. It follows that any
class $\alpha \in \TF_2(\bar{V}|\bar{K};p,\Zp)$ whose image by the map
$$\hat{\Gamma} \colon \TF(\bar{V}|\bar{K};p)\f \to
\tate(\mathbb{T},T(\bar{V}|\bar{K}))\f$$
represents $t^{-1} \in \hat{E}_{2,0}^2$ is a
$\TF_0(\bar{V}|\bar{K};p,\Zp)$-module generator. Let
$$\alpha \in \lim_F T_p W_n\Omega_{(\bar{V},\bar{M})}^1
\xto{\sim} \TF_2(\bar{V}|\bar{K};p,\Zp)$$
to be the unique class such that $\hat{\Gamma}(\alpha) =
\hat{\tau}_{\e}(\lambda_{-L})$. Since $\hat{\beta} = ([L] - 1) \cdot
\lambda_{-L}$ by Prop.~\ref{topologicalktheory} and since
Lemmas~\ref{tauhatbeta} and~\ref{tauhatL} show that the composite
$$\pi_2(\tate(\mathbb{T},f^*ku)) \xto{\hat{\tau}_{\e}}
\pi_2(\tate(\mathbb{T},T(\bar{V}|\bar{K})),\Zp) \xleftarrow{\sim}
\TF_2(\bar{V}|\bar{K};p,\Zp)$$
takes $([L] - 1) \cdot \lambda_{-L}$ to $\psi([\e_1'] - 1) \cdot
\alpha$ and $\hat{\beta}$ to $\beta_{\e} = \beta_{\e}^{\TF}$, we find
that
$$\psi([\e_1'] - 1) \cdot \alpha = \beta_{\e}.$$
We claim that $u = ([\e_1'] - 1)/([\e_1] - 1)$ is a unit in
$W(R_{\bar{V}})$. Indeed, this is the case if and only if the image
$\bar{u} = (\e_1' - 1)/(\e_1 - 1)$ by the canonical projection
$W(R_{\bar{V}}) \to R_{\bar{V}}$ is a unit in $R_{\bar{V}}$. But
$R_{\bar{V}}$ is a valuation ring and
$$v_R(\bar{u}) = v_R(\e_1' - 1) - v_R(\e - 1) = \frac{1}{p-1} -
\frac{1}{p-1} = 0,$$
so $\bar{u}$ is a unit. Hence $\alpha_{\e} = \psi(u) \cdot \alpha$ is
a generator of $\TF_2(\bar{V}|\bar{K};p,\Zp)$ and the formula
$\beta_{\e} = \psi([\e_1] - 1) \cdot \alpha_{\e}$ holds as desired.
\end{proof}

\section{Galois invariants}\label{galoissection}

\subsection{}We show in Prop.~\ref{vanishingcohomology} below that in
positive degrees, the rational homotopy groups
$\TR_q^n(\bar{V}|\bar{K};p,\Qp)$ have vanishing Galois cohomology.

\begin{lemma}\label{invertp}For every ring $A$, the ghost map induces
an isomorphism
$$w \otimes 1 \colon W_n(A) \otimes \Z[\textstyle{\frac{1}{p}}]
\xto{\sim} A^n \otimes \Z[\textstyle{\frac{1}{p}}].$$
\end{lemma}

\begin{proof}Suppose first that $A$ is $p$-torsion free ring, possibly
without unit. Then the ghost map restricts to an isomorphism
$$w\:W_n(pA)\xto{\sim}\prod_{s=0}^{n-1}p^{s+1}A.$$
Hence we have a map of short-exact sequences
$$\xymatrix{
{0} \ar[r] &
{W_n(pA)} \ar[r] \ar[d]^(.4){\sim} &
{W_n(A)} \ar[r] \ar[d] &
{W_n(A/pA)} \ar[r] \ar[d] &
{0} \cr
{0} \ar[r] &
{\displaystyle\prod_{s=0}^{n-1}p^{s+1}A} \ar[r] &
{\displaystyle\prod_{s=0}^{n-1}A} \ar[r] &
{\displaystyle\prod_{s=0}^{n-1}A/p^{s+1}A} \ar[r] &
{0.}
}$$
Since $\Z[\frac{1}{p}]$ is flat over $\Z$, we get an induced map of
short-exact sequences
$$\xymatrix{
{0} \ar[r] &
{W_n(pA)\otimes\Z[\frac{1}{p}]} \ar[r] \ar[d]^{\sim} &
{W_n(A)\otimes\Z[\frac{1}{p}]} \ar[r] \ar[d] &
{W_n(A/pA)\otimes\Z[\frac{1}{p}]} \ar[r] \ar[d] &
{0} \cr
{0} \ar[r] &
{{\displaystyle\prod_{s=0}^{n-1}}p^{s+1}A\otimes\Z[\frac{1}{p}]} \ar[r] &
{{\displaystyle\prod_{s=0}^{n-1}}A\otimes\Z[\frac{1}{p}]} \ar[r] &
{{\displaystyle\prod_{s=0}^{n-1}}A/p^{s+1}A\otimes\Z[\frac{1}{p}]} \ar[r] &
{0.}
}$$
The right hand terms are both zero. Hence, the middle vertical map is
an isomorphism. Finally, for a general ring $A$, we write $A=P/I$
where $P$, and hence $I$, is a ring without $p$-torsion. Then the
diagram
$$\xymatrix{
{0} \ar[r] &
{W_n(I)\otimes\Z[\frac{1}{p}]} \ar[r] \ar[d]^{\sim} &
{W_n(P)\otimes\Z[\frac{1}{p}]} \ar[r] \ar[d]^{\sim} &
{W_n(A)\otimes\Z[\frac{1}{p}]} \ar[r] \ar[d] &
{0} \cr
{0} \ar[r] &
{{\displaystyle\prod_{s=0}^{n-1}}I\otimes\Z[\frac{1}{p}]} \ar[r] &
{{\displaystyle\prod_{s=0}^{n-1}}P\otimes\Z[\frac{1}{p}]} \ar[r] &
{{\displaystyle\prod_{s=0}^{n-1}}A\otimes\Z[\frac{1}{p}]} \ar[r] &
{0}
}$$
completes the proof.
\end{proof}

\begin{lemma}\label{rational}For all $n,q\geq 1$, the canonical
inclusion
$$W_n(\bar V)\f\cdot\beta_{\e,n}^q
\to W_n(\bar V)\f\cdot\alpha_{\e,n}^q$$
becomes an isomorphism after inverting $p$.
\end{lemma}

\begin{proof}By Lemma~\ref{invertp}, the ghost map induces an
isomorphism
$$w\:W_n(\bar V)\f[{\textstyle{\frac{1}{p}}}]\xto{\sim}
\prod_{i=0}^{n-1}\bar K\f.$$
And since $\beta_{\e,n}=\theta_n([\e_n]-1) \cdot \alpha_{\e,n}$, it
suffices to show that ghost map takes the element $\theta_n([\e_n] -
1) = [\e_1^{(n)}]_n - 1$ to a unit. But
$$w([\e_1^{(n)}]_n-1)=([\e_1^{(n)}]_n-1,[\e_1^{(n)}]_n^p-1,\dots,[\e_1^{(n)}]_n^{p^{n-1}}-1),$$
and since $\e_1^{(n)}$ is a primitive $p^n$th root of unity, the ghost
coordinates are all non-zero.
\end{proof}

\begin{prop}\label{vanishingcohomology}Let $q$ be a positive
integer. Then the continuous cohomology group
$H_{\operatorname{cont}}^i(G_K,\TR_q^n(\bar V|\bar K;p,\Qp))$ is zero,
for all $i \geq 0$.
\end{prop}

\begin{proof}By Lemma~\ref{rational}, we have a canonical isomorphism
of $G_K$-modules
$$W_n(\bar V)\f \otimes_{\Zp} \Qp(m) \xto{\sim}
\TR_{2m}^n(\bar{V}|\bar{K};p,\Qp).$$
The $V$-filtration of the left-hand side is finite of length $n$, and
$$\gr_V^sW_n(\bar V)\f\otimes_{\Zp}\Qp(m) = \bar{K}\f(m).$$
But Tate~\cite{tate} has shown that
$H_{\operatorname{cont}}^i(G_K,\bar K\f(m))$ is zero, for all $i \geq
0$, provided that $m > 0$. 
\end{proof}

\providecommand{\bysame}{\leavevmode\hbox to3em{\hrulefill}\thinspace}
\providecommand{\MR}{\relax\ifhmode\unskip\space\fi MR }
\providecommand{\MRhref}[2]{%
  \href{http://www.ams.org/mathscinet-getitem?mr=#1}{#2}
}
\providecommand{\href}[2]{#2}

\end{document}